\documentclass[11pt]{article}

\def\mathbb{\Bbb}
\usepackage{amsfonts,latexsym,amstext}
\usepackage{amsmath}

\usepackage{amssymb}

\newtheorem{theorem}{Theorem}[section]
\newtheorem{lemma}[theorem]{Lemma}
\newtheorem{proposition}[theorem]{Proposition}
\newtheorem{definition}{Definition}[section]

\newtheorem{hypothesis}[theorem]{Hypothesis}

\newtheorem{remark}[theorem]{Remark}
\newtheorem{corollary}[theorem]{Corollary}
\makeatletter
\@addtoreset{equation}{section}

\makeatother

\def\qed{{\hfill\hbox{\enspace${ \square}$}} \smallskip}
\def\sqr#1#2{{\vcenter{\vbox{\hrule height .#2pt \hbox{\vrule
 width .#2pt height#1pt \kern#1pt \vrule
width .#2pt} \hrule height .#2pt}}}}
\def\square{\mathchoice\sqr54\sqr54\sqr{4.1}3\sqr{3.5}3}

\def\ds{\begin{displaystyle}}
\def\eds{\end{displaystyle}}
\def\dis{\displaystyle }
\def\<{\langle }
\def\>{\rangle }

\def\R{\mathbb R}

\def\E{\mathbb E}
\def\P{\mathbb P}

\def\L{\mathbb L}
\def\D{\mathbb D}

\def\cala{{\cal A}}
\def\calb{{\cal B}}

\def\calf{{\cal F}}
\def\calg{{\cal G}}

\def\calt{{\cal T}}
\def\calu{{\cal U}}

\def\call{{\cal L}}
\def\cals{{\cal S}}
\def\bfC{{\bf C}}

\title{Stochastic equations with delay: optimal control via BSDEs and
regular solutions of Hamilton-Jacobi-Bellman equations}

\date{}
\author{Marco Fuhrman\\
Politecnico di Milano,
Dipartimento di Matematica\\
piazza Leonardo da Vinci 32, 20133 Milano, Italy\\
e-mail: marco.fuhrman@polimi.it\\
\\
Federica Masiero, Gianmario Tessitore\\
Dipartimento di Matematica e Applicazioni, Universit\`a di Milano Bicocca\\
via Cozzi 53, 20125 Milano, Italy\\
e-mail: federica.masiero@unimib.it,
gianmario.tessitore@unimib.it}

\begin{document}

\maketitle

\begin{abstract}
We consider an Ito stochastic differential equation with delay, driven by brownian motion, whose
solution,
by an appropriate reformulation,  defines a Markov process $X$ with values in
a space of continuous functions $\bfC$, with generator $\call$. We then consider a backward  stochastic
differential equation depending on $X$, with unknown processes $(Y,Z)$, and we  study properties  of the resulting
system, in particular we identify the process $Z$ as a deterministic functional of $X$.
We next prove that the forward-backward system provides a suitable solution to a class of
parabolic partial differential equations on the space $\bfC$ driven by $\call$, and we apply this result
to prove a characterization of the fair price and the hedging strategy for a financial market
with memory effects.
We also include applications to optimal stochastic control of differential equation with delay:
in particular we characterize optimal controls as feedback laws in terms  the process $X$.
\end{abstract}

\section{Introduction}

We will consider stochastic differential equations with delay (SDDEs for short)
on a finite interval of the form
\begin{equation}\label{eqscalareintro}
\left\{\begin{array}{l} dy_t  = b(t ,y_{t +\cdot})\; dt
+\sigma(t ,y_{t +\cdot})\; dW_t ,\qquad t  \in [0,T],\\
y_\theta=x(\theta), \text{   } \theta \in [-r,0],
\end{array}\right.
\end{equation}
 for an unknown process
$(y_t)_{t\in [-r,T]}$ in $\R^n$. Here $r>0$ is called the delay and
we use the notation $y_{t+\cdot}=(
y_{t+\theta})_{\theta \in [-r,0]}$.
It is customary, and convenient, to introduce the space
$
\bfC=C([-r,0];\R^n)
$ and the
$\bfC$-valued process $X=(X_{t })_{t \in[0,T]}$  defined by
$$
    X_t (\theta)=y_{t +\theta},\qquad
\theta \in[-r,0]
$$
With this notation,
$b(t,\cdot)$ and $\sigma(t,\cdot)$ are functions defined on
$\bfC $ and the equation can be written
$$
    \left\{\begin{array}{l}\dis dy_t =
b(t,{X}_t)\; dt
+\sigma(t,{X}t)\; dW_t,\quad t\in[0,T],
\\\dis
X_t=x\in \bfC.
\end{array}\right.
$$
 SDDEs are a classical subject: in the standard reference  book
\cite{Moh1} (see also \cite{Moh2}) basic results are established: existence and
uniqueness of solutions,
regular dependence on parameters, Markov property of $X$ as a $\bfC$-valued process,
characterization of its generator. In \cite{daza2} long time asymptotics
is studied in detail.

In this paper we will present new results on optimal control problems for SDDEs.
Moreover, since
the Markov character of solutions allows for application of dynamic programming
arguments, we will also prove new results on the corresponding
Hamilton-Jacobi-Bellman equation.
More generally, we will consider a class of semilinear versions of the parabolic Kolmogorov equation
associated to the process $X$. This class includes as a very special case some infinite-dimensional
variants of the Black-Scholes equation for the fair price of an option, of great
interest in mathematical
finance and already considered in \cite{chyo}.

The main tool will be the use of techniques from the theory of backward
stochastic differential equations (BSDEs) in the sense of Pardoux-Peng, first
considered in the nonlinear case in the paper
\cite{PaPe1}. We refer to the monographs \cite{ElMa}, \cite{Pa}  for an exposition
of the basic theory.
The BSDE approach that we follow
consists in addressing equation (\ref{eqscalareintro}), but with
generic initial values $t\in [0,T]$ and $x\in\bfC=C([-r,0];\R^n)$,
and then coupling with another equation
 of backward type, with unknown processes $(Y,Z)$. More precisely
one considers the forward-backward system
\begin{equation}\label{fbsdeintro}
    \left\{\begin{array}{l}\dis dy^{t,x}_\tau =
b(\tau,{X}_\tau^{t,x})\; d\tau
+\sigma(\tau,{X}_\tau^{t,x})\; dW_\tau,\quad \tau\in
[t,T]\subset [0,T],
\\\dis
X_t^{t,x}=x,
\\\dis
 dY_\tau^{t,x}=\psi(\tau,X_\tau^{t,x},Y_\tau^{t,x},Z^{t,x}_\tau)\;d\tau+
 Z^{t,x}_\tau\;dW_\tau,
  \\\dis
  Y_T^{t,x}=\phi(X^{t,x}_T),
\end{array}\right.
\end{equation}
where
$\psi:[0,T]\times\bfC\times\R\times\R^d
\rightarrow\mathbb{R}$ and $\mathbb{\phi}:\bfC\rightarrow\mathbb{R}$
are given functions.
One can then define a (deterministic) function $v:[0,T]\times \bfC\to \R$
setting $v(t,x)=Y_t^{t,x}$ and prove that
\begin{equation}
 \label{identificazioneintro}
Y^{t,x}_{\tau}=
v\left(  \tau,X^{t,x}_{\tau}\right),\qquad
Z^{t,x}_{\tau}=\nabla_0
v\left(  \tau,X^{t,x}_{\tau}\right) \sigma(\tau,X^{t,x}_{\tau}) ,
\end{equation}
where $\nabla_0 $  is a differential operator defined by
\begin{equation}\label{defnablazerointro}
\nabla_0v(t,x)=\nabla_xv(t,x) (\{0\})
\end{equation}
where the gradient $\nabla_x v(t,x)$ at point
$(t,x)\in [0,T]\times \bfC$ is an element of the dual space $\bfC^*$, hence 
 an $n$-tuple  of finite Borel measures on $[-r,0]$. Thus,
 $\nabla_0v(t,x)$ is a vector in $\R^n$ whose components
 are the masses at point $0$ of
the components of $\nabla_xv(t,x)$.

It turns out that $v$ is a solution of a
 semilinear
parabolic equation of the  form:
\begin{equation}\label{nlkdegintro}
  \left\{\begin{array}{l}\dis
\frac{\partial v(t,x)}{\partial t}+\call_t v(t,x) =
\psi (t,x,v(t,x),\nabla_0 v(t,x)\,\sigma(t ,x)),\\ \\
\dis v(T,x)=\phi(x), \qquad t\in [0,T],\;x\in \bfC,
\end{array}\right.
\end{equation}
where
 $\call_t$ is the generator of the Markov process $(X^{t,x}_{\tau})$.
In the finite-dimensional case, this was first proved in \cite{PaPe} for classical
solutions of (\ref{nlkdegintro}).
If one considers the controlled SDDE
\begin{equation}\label{SDDEcontrollataintro}
\left\{
\begin{array}{l}
dy^u_s =
b(s,{X}^u_s)\; ds
+\sigma(s,X^u_s)\,
[h(s,X^u_s,u_s))\;ds+\; dW_s]
,\quad s\in
[{t},T],
\\
{X}_{{t}}={x}.
\end{array}
\right.
\end{equation}
where the solution depends on a control process $u(\cdot)$ taking values
in a space $U$,
and $h:\left[
0,T\right]  \times \bfC \times U\to \R^{d}  $ is given, and
one tries to minimize a
 cost functional
\begin{equation}\label{funzcostofinaleintro}
J(t,x,u(\cdot))
=\E\int_t^Tg( u_s)\; ds
+
\E\,\phi(X^u_T),
\end{equation}
where $g:U\to [0,\infty)$, then equation
(\ref{nlkdegintro}) is the associated Hamilton-Jacobi-Bellman equation,
provided
 the hamiltonian function
$\psi:[0,T]\times \bfC\times\R^d\to\mathbb{R}$ is defined by the formula
$$
\psi\left(t, x,z\right)  =\inf\left\{g(u)
+zh\left(t,x,u\right)  :u\in U\right\},
\qquad t\in [0,T], x\in \bfC,z\in\R^d.
$$
This way we eventually prove that $v$ coincides with the value function of the
control problem.

Although BSDEs were known to be useful tools in the study
of control problems and nonlinear partial differential equations, applications to
infinite-dimensional state spaces are more recent and difficult: see e.g.
\cite{fute}, \cite{fute2} for the case of a Hilbert space, and
\cite{Mas} for some related results on Banach spaces. In these papers, as well as in the
present one, the solution of (\ref{nlkdegintro}) is understood in the so-called mild
sense.
Special difficulties are given by the fact that the state
space $\bfC$ is used as the state space of the basic stochastic process $X$.
The reason for doing this is to allow  for great generality on the coefficients $b,\sigma$
of the SDDEs as well as on the cost functional of the control problem. For instance the
functional $\phi$ occurring in (\ref{funzcostofinaleintro}) could have the form
\begin{equation}\label{esempiononlinearita}
    \phi (z)= \int_{[-r,0]}g(z(\theta))\,\mu(d\theta),
\qquad z\in\bfC,
\end{equation}
for some  $g\in C^1(\R)$
and some measure $\mu$ on $[-r,0]$.
The special case when $\mu$ is supported on a finite number of points
is of particular interest and could be studied by direct methods, but it is included
in our results.
More generally, if
 $\phi_1,\ldots,\phi_n$ are functionals with the form (\ref{esempiononlinearita})
 corresponding to functions $g_1,\ldots,g_n\in C^1(\R)$,
 and if $h\in C^1(\R^n)$, then the functional
$
\phi(z)=h(\phi_1(z),\ldots,\phi_n(z)),
$
can also be treated by our methods.
 One could avoid the
use of the space $\bfC$ by looking at $X$ as a process with values in
the space $L^2([0,T];\R^n)$ instead. This was the approach taken in
\cite{fute3}.
However, this leads to restrictions on the applicability of the corresponding results.

Optimal control problems for SDDEs have been thoroughly investigated in recent
years: we refer to the book \cite{ch} and the references therein.
One of the main results there is the characterization of the value function
as the unique viscosity solution of the
Hamilton-Jacobi-Bellman equation (\ref{nlkdegintro}).
This is achieved for controlled equations of
 more general form than (\ref{SDDEcontrollataintro}), in particular with possibly
 control-dependent diffusion coefficients.
 In our paper we assume stronger conditions, namely a special form for the control system
 (\ref{SDDEcontrollataintro}) and differentiability assumptions on the
 data $b,\sigma, \phi,\psi$ with respect
 to the space variable $x\in \bfC$. In this case we are able to prove further properties on
 the solution $v$, namely its  differentiability.
 Thus, some of our results can be viewed as regularity properties,
occurring under appropriate assumptions, of the viscosity solutions
 of (\ref{nlkdegintro}). However, the existence of the gradient of $v$ is of special
 interest in optimal control theory, since it allows to characterize optimal
 controls via feedback laws and to prove existence of optimal controls
 after appropriate formulation.

Parabolic equations on the space $\bfC$ of the form (\ref{nlkdegintro}) have also
been considered for other purposes, in particular as a generalization of
the Black-Scholes equation for the fair price of an option, in case the market
models exhibits memory effects, see \cite{chyo}. In particular, the special
operator $\nabla_0$ defined in (\ref{defnablazerointro})
also occurs in the class of equations considered there.
 In the same spirit in
\cite{AHMP} some formulae of Black-Scholes  type are proved.
Here again the approach based on BSDEs gives new results in comparision
to \cite{chyo}.

The plan of the paper is as follows: in section \ref{sez-prelim} we introduce
notation and review some results on SDDEs, adding some precision on regularity
properties of the solution, concerning in particular  their Malliavin derivative.
Section \ref{joint-quadr-var} is devoted to proving
Theorem \ref{esistejointvar}, which is the
key of many subsequent results; here the operator $\nabla_0$ is introduced.
In section \ref{sezioneforback} we present the forward-backward system
(\ref{fbsdeintro}) and prove in particular the second formula in
(\ref{identificazioneintro}). In
section \ref{sez-control} we study the optimal control problem, we prove  in particular
the so-called fundamental relation, we give criteria for optimality
of feedback controls and we prove existence of optimal controls
in the weak sense.
Section \ref{sez-nlk} is devoted to the study of equation
 (\ref{nlkdegintro}): it is proved that a solution exists and it
 is connected to the solution of the forward-backward
 system (\ref{fbsdeintro}) by formula (\ref{identificazioneintro});
 in particular it is proved that the value function
of the control problem is a solution to this equation  (in the mild sense);
finally, it is shown how  (\ref{nlkdegintro}) may arise as the a
Black-Scholes  equation in a financial market and we give explicit
conditions for its solvability.

\section{Preliminary  results on stochastic delay differential equations}
\label{sez-prelim}
\subsection{Notations}

In this paper we consider a  complete probability space
$\left(  \Omega,\mathcal{F}, \mathbb{P}\right)  $ and a standard
 Wiener process  $W=(W_t)_{t\ge 0}$ with values in $\R^d$.
 We denote by
 $(\calf_t)_{t\ge 0}$ the natural filtration of
 $W$  augmented in the usual way by the sets of
$\P$-measure $0$.

 For fixed $r>0$, we introduce the space
$$
\bfC=C([-r,0];\R^n)
$$
of continuous functions from $[-r,0]$ to $\R^n$,
endowed with the usual norm $|
f|_{\bfC}=\sup_{\theta\in[-r,0]}\vert f(\theta)\vert$. We will consider
$\bfC$-valued stochastic processes: for $T>0$ we
 say that a  $\bfC$-valued process $(X_t)_{t\in[0,T]}$
belongs to the space $\cals ^p([0,T];\bfC)$ ($1\le p<\infty$) if
its path are $\bfC$-continuous $\P$-a.s. and the norm
\[
\|X\|^p_{\cals ^p([0,T];\bfC)}=
 \E \sup_{t\in [0,T]}\vert X_t\vert ^p_{\bfC}=
 \E \sup_{t\in [0,T]}\sup_{\theta \in [-r,0]} \vert X_t(\theta)\vert ^p
\]
is finite. Here and in the following, if
no confusion is possible, we denote the norm of $\R^n$, $\R^d$ and
$\R^{nd}$ by $\vert \cdot \vert$.

We next define several classes of differentiable functions between Banach
spaces, first introduced in \cite{fute} in connection with stochastic processes,
which allow to formulate several regularity results in a compact way.

In the following, if $E$ and $K$ are Banach spaces, we denote by
$\calg^1(E,K)$ the space of continuous functions $u:E\to K$ such that:
1) $u$ is continuous; 2) $u$ is
G\^ateaux  differentiable on $E$, with G\^ateaux  differential
at point $x\in E$ denoted by $\nabla u(x)\in L(E,K)$ (the latter
being the space of bounded linear operators from $E$ to $K$, endowed
with its usual norm); 3) for every $h\in E$, the map $x\to
\nabla u(x)\,h$ is continuous from $E$ to $K$.
We note that the map $x\to \nabla u(x)$ is not required to be continuous
from $E$ to $L(E,K)$: if this happens then $u$ is also Fr\'echet differentiable.

We say that a function
$v:[0,T]\times E\to K$ belongs to $\calg^{0,1}([0,T]\times E,K)$ if:
1) $v$ is continuous; 2) for every $t\in [0,T]$, $v(t,\cdot)$ is
G\^ateaux  differentiable on $E$, with G\^ateaux  differential
at point $x\in E$ denoted by $\nabla_x v(t,x)\in L(E,K)$;
3) for every $h\in E$, the map $(t,x)\to
\nabla_x v(t,x)\,h$ is continuous from $[0,T]\times E$ to $K$.

Now suppose  $E=C([a,b];\R^n)$, where $a,b\in\R$,
$a<b$.
 We recall that the dual space of $C([a,b])$ is the space of finite Borel measures
 on $[a,b]$, endowed with the variation norm. Identifying $E$ with the product space
  $C([a,b])^n$ in the obvious way we conclude that the dual space $E^*$ of $E$ can be
identified with the space of  $n$-tuples $\mu=(\mu_k)_{k=1}^n$,
where each $\mu_k$ is a finite Borel measure on $[a,b]$, and the value of
$\mu$ at an element $g=(g_k)_{k=1}^n\in C([a,b])^n$, where
$g_k\in C([a,b])$, is denoted
$$
\int_{[a,b]}g(\theta)\cdot \mu(d\theta)=
\sum_{k=1}^n\int_{[a,b]}g_k(\theta)\,\mu_k(d\theta).
$$

Let $v:[0,T]\times \bfC\to \R$ be a  function such
that $v(t,\cdot)$ is G\^ateaux differentiable on $\bfC$ for every $t\in[0,T]$.
Then the gradient $\nabla_x v(t,x)$ at point
$(t,x)\in [0,T]\times \bfC$ is
 an $n$-tuple  of finite Borel measures on $[-r,0]$. We denote by
$|\nabla_x v(t,x)|$ its total variation norm and we
define
\begin{equation}\label{defnablazero}
\nabla_0v(t,x)=\nabla_xv(t,x) (\{0\})
\end{equation}
i.e., $\nabla_0v(t,x)$ is a vector in $\R^n$ whose components
$\nabla_0^kv(t,x)$ ($k=1,\ldots,n$) are the masses at point $0$ of
the components of $\nabla_xv(t,x)$.

\begin{remark}\label{spaziosol}{\em
In the following, a basic role will be played by the space
$\calg^{0,1}([0,T]\times \bfC,\R)$: according to the previous
definitions, it consists of real continuous functions $v$ on $[0,T]\times \bfC$
such that,
for every $t\in [0,T]$, $v(t,\cdot)$ is
G\^ateaux  differentiable on $\bfC$, with G\^ateaux  differential
at point $x\in \bfC$ denoted by $\nabla_x v(t,x)$
(an $n$-tuple of finite Borel measures on $[-r,0]$), such that
 the map
$$(t,x)\to
\<\nabla_xv(t,x),h\>_{\bfC^*,\bfC}= \int_{[-r,0]}h(\theta)\cdot \nabla_xv(t,x)(d\theta)
$$
is continuous on $[0,T]\times \bfC$, for every $h\in \bfC$.
}
\end{remark}

\subsection{Stochastic delay differential equations}

We fix $T>0$ and  we consider the following
stochastic delay differential equation for an unknown process
$(y_t)_{t\in [0,T]}$ taking values in $\R^n$:

\begin{equation}\label{eqscalare}
\left\{\begin{array}{l} dy_t  = b(t ,y_{t +\cdot})\; dt
+\sigma(t ,y_{t +\cdot})\; dW_t ,\qquad t  \in [0,T],\\
y_\theta=x(\theta), \text{   } \theta \in [-r,0],
\end{array}\right.
\end{equation}
where $y_{t+\cdot}$ denotes the past trajectory from time $t-r$
up to time $t$, namely $y_{t+\cdot}=(
y_{t+\theta})_{\theta \in [-r,0]}$, and $r>0$ is the
delay. $b(t,\cdot)$ and $\sigma(t,\cdot)$ are functions of the past
trajectory of $y$ and they are defined on the space of
continuous functions, namely $b:[0,T]\times \bfC
\rightarrow \R^n$ and $\sigma:[0,T]\times \bfC
\rightarrow \R^{nd}$, where $\R^{nd}$ is identified with $L(\R^d,
\R^n)$ the space of linear operators from $\R^{d}$ to $\R^n$.  The function $x\in
\bfC $ is the initial condition. We will refer to
equation (\ref{eqscalare}) as {\em delay equation}.

We make the following assumptions on the coefficients of (\ref{eqscalare}).
\begin{hypothesis}\label{ipotesi1}
\begin{enumerate} \item The functions
$b:[0,T]\times \bfC
\rightarrow \R^n$ and $\sigma:[0,T]\times \bfC
\rightarrow \R^{nd}$ are continuous
and there exists a constant $K>0$ such that
for all
$t\in [0,T]$ and
 $y(\cdot)\in \bfC$
$$
|b(t,y(\cdot))|+|\sigma(t,y(\cdot))|\le K\,
(1+|y(\cdot)|_\bfC);
$$
\item there exists a constant $L>0$ such that for all
$t\in [0,T]$ and
 $y(\cdot),z(\cdot)\in \bfC$
$$
|b(t,y(\cdot))-b(t,z(\cdot))|+|\sigma(t,y(\cdot))-
\sigma(t,z(\cdot))|\le L\,
|y(\cdot)-z(\cdot)|_\bfC;
$$
\item for all $t\in[0,T]$, $b(t,\cdot)\in \calg^1(\bfC, \R^n)$,
$\sigma(t,\cdot)\in \calg^1(\bfC, \R^{nd})$.
\end{enumerate}
\end{hypothesis}

In the following we collect
some results on existence and
uniqueness of a solution to equation (\ref{eqscalare}) and on its
regular dependence on the initial condition.
It turns out that there exists a continuous solution, so we can define
a $\bfC$-valued process $X=(X_{t })_{t \in[0,T]}$  by
\begin{equation}\label{defdix}
    X_t (\theta)=y_{t +\theta},\qquad
\theta \in[-r,0].
\end{equation}
We notice that if $t +\theta <0$
then $y_{t +\theta}=x(t +\theta)$. We will use the notations $y^x$ , $y^x_{t }$,
$X^x$ or $X^x_{t }$ to indicate dependence on the
starting point $x\in \bfC$.

\begin{theorem}\label{teoSDDE}
If Hypothesis \ref{ipotesi1}, points 1 and 2, holds true then there
exists a unique continuous adapted solution
of the delay equation (\ref{eqscalare}), and moreover the process  $(X_t)_{t\in[0,T]}$
belongs to
$ \cals ^p([0,T];\bfC)$
for every $p\geq 2$ and
\[
\|X\|^p_{\cals ^p([0,T];\bfC)}=
 \E \sup_{t\in [-r,T]}\vert y_{t}\vert ^p \leq C
\]
for some constant $C>0$ depending only on $K,L,T,p$.

In addition,
the map $x\to X^x$
is Lipschitz continuous  from $\bfC$ to
$\cals^p([0,T];\bfC)$; more precisely,
\[
\Vert X^{x_1}-X^{x_2}\Vert_{\cals ^p([0,T];\bfC)}
=  \left(\E \sup_{t\in [-r,T]}\vert y^{x_1}_t-y^{x_2}_{t}\vert ^p \right)^{1/p}
\leq L\sup_{\theta\in [-r,0]}\vert x_1(\theta)-x_2(\theta)\vert
\]
for some constant $L>0$ depending only on $K,L,T,p$.

If we further assume that Hypothesis \ref{ipotesi1}, point 3,
holds true then the map $x\to X^x$ belongs to the space
 $\calg^1 (\bfC, \cals^p([0,T]; \bfC))$.

\end{theorem}

\noindent {\bf Proof.} For the proof (in the case of $p=2$), we
refer to \cite{Moh1}, Chapter II:  we refer to Theorem 2.1
for existence and uniqueness of the solution of equation
(\ref{eqscalare}), to Theorem 3.1 for the Lipschitz dependence of
this solution on the initial datum, and to Theorem 3.2 for the
differentiability of the solution with respect to the initial
datum. See also \cite{Moh2}, Theorems I.1 and  I.2. The
proof in the case of $p>2$ can be performed in a similar way. \qed

Let us introduce a delay equation similar to (\ref{eqscalare}) but with initial
condition given at time $t\in[0,T]$:

\begin{equation}\label{eqscalare-s}
\left\{\begin{array}{l} dy_{\tau}^{t,x} =
b(\tau,y_{\tau+\cdot}^{t,x}) \;d{\tau}
+\sigma(\tau,y_{\tau+\cdot}^{t,x})\; dW_\tau,\qquad \tau \in [t,T],\\
y_{t+\theta}^{t,x}=x(\theta), \qquad \theta \in [-r,0].
\end{array}\right.
\end{equation}
We introduce  the $\bfC$-valued process given by
\begin{equation}\label{defdixdue}
X_\tau^{t,x}(\theta)=y_{\tau+\theta}^{t,x}, \qquad \theta\in
[-r,0].
\end{equation}
By \cite{Moh1}, Chapter III, Theorem 2.1, the $\bfC$-valued process
$(X_{\tau}^{t,x})_{\tau\in[t,T]}$ is a Markov process with
transition semigroup, acting on bounded and Borel measurable  $\phi:\bfC\to \R$,
given by
\begin{equation}\label{defdiP}
P_{t, \tau}[\phi](x)=\E\, \phi(X_{\tau}^{t,x}),
\qquad 0\le t\le \tau\le T,\;x\in\bfC.
\end{equation}

\begin{remark}\label{remarkL}{\em
The transition semigroup $(P_{t, \tau})$ has been extensively studied in the literature,
see e.g. \cite{Moh1} and \cite{Moh2}. For the sake of completeness,
we briefly recall some result on its generator, which will appear in section 6
in the formulation of the Kolmogorov equation.
For simplicity, let us consider the autonomous case in
equation (\ref{eqscalare-s}): $b$ and $\sigma$ do not depend on
time and $s=0$, so we consider the one parameter semigroup
$(P_t)_{t\in [0,T]}$. The transition semigroup
$(P_t)$ is never strongly continuous on the space $\bfC$,
nevertheless it admits a weakly continuous generator $\call$, see
\cite{Moh1},
 chapter IV and \cite{Moh2}, chapter II. Let
$S_t:\bfC\rightarrow\bfC$ denote  the shift operator and
let $\cals$ denote the weak generator of the corresponding
 semigroup. To derive a formula for the generator $\call$ we need to
augment $\bfC$ by adding an $n$-dimensional direction.
$\call$ will be equal to the sum of the generator of the shift semigroup
 $\cals$ and a second order linear partial differential
operator along this new direction. Let $F_n:=\left\lbrace
v1_{{0}}:v\in\R^n\right\rbrace $ and $\bfC\oplus F_n:=\left\lbrace
f+v1_{{0}}:f\in \bfC, v\in\R^n\right\rbrace $ with the norm $\Vert
f+v1_{{0}} \Vert _{\bfC\oplus F_n}:=\vert f\vert_{\bfC}+\vert
v\vert$. Suppose that $\phi:\bfC\rightarrow \R$ is twice
continuously Frechet differentiable and let $f\in \bfC$. Then the
Frechet derivatives $\nabla \phi(f)$ and  $\nabla^2 \phi(f)$ have
unique weakly continuous linear and bilinear extensions
\begin{equation*}
 \overline{ \nabla \phi(f)}:\bfC\oplus F_n\rightarrow \R,\qquad
 \overline{\nabla^2 \phi(f)}:(\bfC\oplus F_n)\times(\bfC\oplus F_n)\rightarrow \R.
\end{equation*}
We are ready to introduce $\call$. Suppose that $\phi
:\bfC\rightarrow \R$, $\phi\in D(\cals)$, and $\phi$ is
sufficiently smooth (e.g. $\phi$ is twice continuously
differentiable and its derivatives are globally bounded and
lipschitz continuous). Then $\phi\in D(\call)$ and $\forall f \in
\bfC$
\begin{equation}\label{weakgenerator}
 \call (\phi)(f)=\cals(\phi)(f)+\overline{ \nabla \phi(f)}(b(f)1_{{0}})
+\dfrac{1}{2}\sum_{i=1}^n \overline{\nabla^2 \phi(f)}(\sigma(f)(e_i)1_{{0}},\sigma(f)(e_i)1_{{0}}),
\end{equation}
where $\left\lbrace e_i\right\rbrace
_{i=1}^n$ is any basis of $\R^n$.
}\end{remark}

\subsection{Differentiability in the Malliavin sense}

Our  aim is now to compute the Malliavin derivative
of the solution of the delay equation.
We start by recalling  some basic definitions from the Malliavin calculus.
We refer the reader to the book \cite{Nu} for a detailed
exposition.

We consider again a  standard
 Wiener process  $W=(W_t)_{t\ge 0}$ in $\R^d$ and
the Hilbert space $L^2([0,T];\R^d)$ of Borel
measurable, square summable functions on $[0,T]$ with values in
$\R^d$, with its natural inner product. This can be identified
with the product space $(L^2([0,T]))^d$ or with the space
$L^2(\calt)$, where the measure space
$\calt:=[0,T]\times\{1,\ldots,d\}$ is endowed with the product of
the Lebesgue measure on $[0,T]$ and the counting measure on
$\{1,\ldots,d\}$. Elements $h\in L^2([0,T];\R^d)$ may be denoted
$\{h^j(s),\; s\in [0,T],\,j=1,\ldots,d\}$ or
 $\{h^j\}$, where $h^j\in L^2([0,T])$.

For every $h\in L^2([0,T];\R^d)$ we denote
$$W(h)=\int_0^T h(s)\cdot\,dW_s=
\sum_{j=1}^dh^j(s)\,dW^j_s.
$$
$W$ is an isometry of $L^2([0,T];\R^d)$ onto a gaussian subspace
of $L^2(\Omega)$, called the first Wiener chaos. Given a Hilbert
space $K$, let $S_K$ be the set of $K$-valued random variables $F$
of the form
$$
F=\sum_{r=1}^mf_r(W(h_1),\ldots,W(h_n))\, k_r,
$$
where
 $h_1,\ldots , h_n\in L^2([0,T];\R^d)$,
  $\{k_r\}$ is a basis of $K$ and $f_1,\ldots f_m$ are
infinitely differentiable functions $\R^n\to \R$ bounded together
with all their derivatives. The Malliavin derivative $DF$ of $F\in
S_K$ is defined as the process $\{D_s^jF;\; s\in [0,T],\,
j\in\{1,\ldots,d\}\}$ given by
$$D^j_sF= \sum_{r=1}^m\sum_{k=1}^n\partial_k f_r
(W(h_1),\ldots,W(h_n))\,  h_k^j(s)\,k_r,
$$
with values in $K$; by $\partial_k$ we denote the partial
derivatives with respect to the $k$-th variable. It is known that
the operator $D:S_K\subset L^2(\Omega;K)\to L^2(\Omega\times [0,T]
\times \{1,\ldots,d\};K)=L^2(\Omega\times \calt ;K)$ is closable.
We denote by $\D^{1,2}(K)$ the domain of its closure, endowed with
the graph norm, and we use the same letter to denote $D$ and its
closure:
$$
D:\D^{1,2}(K) \subset L^2(\Omega;K) \to L^2(\Omega\times \calt;
K).
$$
The adjoint operator of $D$,
$$
\delta: {\;\rm dom\;}(\delta)\subset L^2(\Omega\times \calt;K)\to
L^2(\Omega;K),
$$
is called Skorohod integral. For a process $u=\{u_s^j;\; s\in
[0,T],\, j\in\{1,\ldots,d\}\} \in$ dom$(\delta)$ we will also use
the notations
$$
\delta (u)=\int_0^Tu_s\;\hat{d}W_s=
\sum_{j=1}^d\int_0^Tu_s^j\;\hat{d}W^j_s.
$$
It is known that dom$(\delta)$ contains every
$(\calf_t)$-predictable process in $L^2(\Omega\times \calt;K)$ and
for such processes the Skorohod integral coincides with the It\^o
integral; dom$(\delta)$ also contains  the class $\L^{1,2}(K)$,
the latter being defined as the space of processes $u\in
L^2(\Omega\times \calt;K)$ such that $u_t^j\in \D^{1,2}(K)$ for
a.e. $t\in [0,T]$ and every $j$,
 and there exists a
measurable version of $D_s^iu_t^j$ satisfying
$$
\| u \|^2_{\L^{1,2}(K)}=
 \| u \|^2_{L^2(\Omega\times \calt;K)}
+\E\sum_{i,j=1}^d\int_0^T \int_0^T \|D_s^i
u^j_t\|^2_{K}\,dt\,ds<\infty.
$$
Moreover, $\| \delta(u) \|^2_{ L^2(\Omega;K)}\leq \| u
\|^2_{\L^{1,2}(K)}$. We note that the space $\L^{1,2}(K)$ is
isometrically isomorphic to $L^2(\calt ; \D^{1,2}(K))$.

Finally we recall that if $F\in \D^{1,2}(K)$ is measurable with
respect to $\calf_t$ then  $D^jF=0$ a.s. on $\Omega \times (t,T] $
for every $j$.

 If $K=\R$ or $K=\R^n$, we write $\D^{1,2}$ and $\L^{1,2}$
instead of $\D^{1,2}(K)$ and $\L^{1,2}(K)$ respectively.

We now introduce the Malliavin derivative for a functional of a stochastic process.
In the remainder of this section we set $E=C([-r,T];\R^n)$.
If $f\in \calg^1(E, \R^n)$ then, according to the notation introduced above,
$$
\< \nabla f(x), g\>_{E^*,E}=\dis \int_{[-r,T]}g(\theta)\cdot \nabla f(x)(d\theta),
\qquad x,g\in E.
$$
If $y$ is a continuous stochastic process with time parameter $[-r,T]$ then
$f(y_.)$ is a random variable. We wish to state  a chain rule
for the Malliavin derivative of $f(y_.)$. We will restrict to the case
when $y$ is adapted, more precisely its restriction to $[0,T]$ is
adapted to $(\calf_t)_{t\in [0,T]}$  and its restriction to $[-r,0]$ is
{\em deterministic}. Clearly, $Dy_t=0$ for $t\in[-r,0]$.
Following \cite{Hu}, lemma 2.6, we have the following basic
result (we note that in \cite{Hu} derivatives are understood in the
sense of Fr\'echet, but the same arguments apply to the present
situation).

\begin{lemma}\label{chainrule}
For
$E=C([-r,T];\R^n)$, let
$f\in \calg^1(E, \R)$ be a Lipschitz continuous function. Assume that  $y=(y_t)_{t\in [-r,T]}$
is a process in $\R^n$ satisfying the following conditions:
\begin{enumerate}
\item $y$ is a continuous
adapted process and
$
\E \sup _{t\in [-r,T]}\vert y_t\vert^2<\infty$;
\item
$y\in  L^2([-r,T],\D^{1,2})$ and  the process
$\{D_sy_t, 0\le s\le t\le T\}$ admits a version such that,
for every $ s\in [0,T]$, $\{D_sy_t, t\in [s,T]\}$ is a continuous process and
$$
\E \dis \int_{0}^T\sup _{t\in [s,T]}\vert D_sy_t\vert^2ds<\infty.
$$
\end{enumerate}
Then $f(y_.)\in \D^{1,2}$ and its Malliavin derivative is given by the formula:
for  $j=1,...,d$ and a.e. $s\in [0,T]$ we have, $\P$-a.s.,
\begin{equation}\label{formulachainrule}
 D^j_s(f(y_{\cdot}))=\<\nabla f(y_.),D^j_sy_{\cdot}\>_{ E^*,E}=\int_{[-r,T]}D^j_sy_\theta\cdot
 \nabla f(y_.)(d\theta).
\end{equation}

\end{lemma}

Next we establish when the solution of the delay equation is Malliavin differentiable,
and moreover we write a stochastic (functional) differential equation satisfied
by the Malliavin derivative. We substantially follow \cite{Hu}, Theorem 4.1.
\begin{theorem}\label{derdimalliavinSDDE}
 Let Hypothesis
\ref{ipotesi1} be satisfied. Then the solution $(y_t)_{t\in [-r,
T]}$ satisfies conditions 1. and 2. in
Lemma \ref{chainrule}. Moreover $y_t\in \D^{1,2}$ for every $t\in [0,T]$ and
the following equation holds:
for  $j=1,...,d$ and every $s\in [0,T]$ we have, $\P$-a.s.,
\begin{equation}\label{eqdermallavinSDDE}
\left\lbrace \begin{array}{lll} D^j_sy_t
&=&\dis \sigma(s,y_{s+\cdot})
+\dis\int_s^t\int_{[-r,0]}D^j_sy_{t+\theta}\cdot\nabla_x b(t,y_{t+\cdot})(d\theta)\; dt\\
&&+\dis\int_s^t\int_{[-r,0]}D^j_sy_{t+\theta}\cdot\nabla_x \sigma(t,y_{t+\cdot})(d\theta)\; dW_t,
\qquad t \in [s,T],\\
D^j_sy_t&=&0,\qquad t \in [-r,s).
\end{array}\right.
\end{equation}
Finally, for every $p\in [2,\infty)$ and $s\in [0,T]$ we have
\begin{equation}\label{stimalpderivmall}
\E \dis \int_{0}^T\sup _{t\in [s,T]}\vert D_sy_t\vert^pds<\infty.
\end{equation}
\end{theorem}
\noindent {\bf Proof.} Except for
the final statement,
the proof can be achieved with techniques similar to the ones
indicated in the proof of Theorem 4.1 in \cite{Hu}. The only minor  difference
is that we consider a general delay differential equation, while in \cite{Hu}
the coefficients depend on the past behavior of the
solution only after time $0$: however, the same arguments apply.

The proof of the final statement follows by standard estimates
on equation (\ref{stimalpderivmall}), taking into account that
$\nabla_xb$ and $\nabla_x\sigma$ are bounded in the total variation norm.
\qed

\begin{corollary}\label{chainruleperlasol}
Suppose that the assumptions of
Theorem \ref{derdimalliavinSDDE} hold true, let $\bfC=C([-r,0];\R^n)$ and
$(X_t)_{t \in[0,T]}$ be the $\bfC$-valued process defined by (\ref{defdix}).
Suppose that  $f\in \calg^1(\bfC;
\R)$ satisfies
$$
   |\nabla f(x)|\leq C(1+|x|_\bfC)^{m}, \qquad x\in \bfC,
   $$
   for some $C>0$ and $m\ge 0$.

Then for every $t\in [0,T]$,
$f(X_t)=f(y_{t+\cdot})$ belongs to $\D^{1,2}$ and for
$j=1,...,d$ we have, for
a.e. $s\in [0,T]$, $\P$-a.s.,
\begin{equation}\label{formulachainruleperlasol}
 D^j_s(f(X_{t}))=\<\nabla f(X_{t}),D^j_sy_{t+\cdot}\>_{\bfC^*,\bfC}
 =\int_{[-r,0]}D^j_sy_{t+\theta}\cdot\nabla f(X_{t})(d\theta).
\end{equation}
\end{corollary}

\noindent {\bf Proof.} The conclusion follows immediately from Lemma \ref{chainrule}
and Theorem \ref{derdimalliavinSDDE}
if $f$ is a Lipschitz function. The general case can be proved by approximating $f$
by a sequence of Lipschitz functions obtained by a standard truncation procedure.
\qed

\begin{remark}{\em
The first result on Malliavin differentiability of the solution of a functional
stochastic differential equations was proved in \cite{KS}.
In that paper  the aim was to prove that $y_t$ belongs to the domain
of the generator
of the Ornstein-Uhlenbeck semigroup of the Malliavin calculus, therefore
more restrictive assumptions were  assumed on the
 coefficients of equation (\ref{eqscalare}), in particular they
were required to be twice differentiable.
}\end{remark}

\section{A result on joint quadratic variations}\label{joint-quadr-var}

The aim of this section is to state and prove
 a  technical result,  Theorem \ref{esistejointvar}, which will be used
in the rest of this paper. To state this theorem we
need to recall some definitions concerning joint quadratic variations
of stochastic processes and to introduce a differential operator,
denoted  $\nabla_0$, which will also play a basic role in the sequel.

We say that a pair of real stochastic processes $(X_t,Y_t)$,
$t\ge0$, admits a joint quadratic variation on the interval
$[0,T]$ if setting
$$
C^\epsilon_{[0,T]}(X,Y)=\frac{1}{\epsilon} \int_0^T
(X_{t+\epsilon}-X_t)(Y_{t+\epsilon}-Y_t)\;dt, \qquad \epsilon>0,
$$
the limit $\lim_{\epsilon\to 0}C^\epsilon_{[0,T]}(X,Y)$ exists in
probability. The limit will be denoted $\<X,Y\>_{[0,T]}$.

This definition is taken from \cite{ruva2}, except that we do not
require that the  convergence in probability holds uniformly with
respect to time. In \cite{ruva2} the process $\<X,Y\>$ is called
generalized covariation process; several properties are
investigated in \cite{ruva3}, \cite{ruva4}, often in connection
with the stochastic calculus introduced in \cite{ruva1}. With
respect to the classical definition, the present one has some
technical advantages that are useful when dealing with convergence
issues (compare for instance the proof of Theorem
\ref{esistejointvar} below).

In the following we will consider joint quadratic variations over
different intervals, which is defined by obvious modifications.

It is easy to show that if $X$  has paths with finite variation
and  $Y$ has continuous paths then $\<X,Y\>_{[0,T]}=0$.

If $X$ and $Y$ are stochastic integrals with respect to the Wiener
process then the joint quadratic variation as defined above
coincides with the classical one. A similar conclusion holds for
general semimartingales:  see \cite{ruva2}, Proposition 1.1.

We set $\bfC=C([-r,0];\R^n)$ and, for every $t\in [0,T]$ and $x\in
\bfC$, we let $\{X_s^{t,x},\; s\in [t,T]\}$ denote the process
defined by the equality (\ref{defdixdue}), obtained as a solution to
equation (\ref{eqscalare-s}).
In particular it is an $\bfC$-valued process with continuous paths
and adapted to the filtration $\{\calf_{[t,s]},\; s\in [t,T]\}$.
$X_s^{t,x}(\omega)$ is measurable in $(\omega, s,t,x)$.

Let $u:[0,T]\times \bfC\to \R$ be a  function such
that $u(t,\cdot)$ is G\^ateaux differentiable on $\bfC$ for every $t\in[0,T]$.
Then the gradient $\nabla_x u(t,x)$ at point
$(t,x)\in [0,T]\times \bfC$ is
an $n$-tuple  of finite Borel measures on $[-r,0]$;  we denote by
$|\nabla_x u(t,x)|$ its total variation norm and we denote
$
\nabla_0u(t,x)=\nabla_xu(t,x) (\{0\})
$, compare (\ref{defnablazero}); thus, $\nabla_0u(t,x)$ is a vector in $\R^n$ whose components
$\nabla_0^ku(t,x)$ ($k=1,\ldots,n$) are the masses at point $0$ of
the components of $\nabla_xu(t,x)$.

We denote $W^i$ ($i=1,\ldots,d$) the $i$-th component of the
Wiener process $W$,
 by $\sigma^i$ the $i$-th column of the $n\times d$ matrix
$\sigma$, and by $\sigma^i_k$ ($k=1,\ldots, n$) its components.

\begin{theorem}\label{esistejointvar}
  Assume that  $u:[0,T]\times \bfC\to \R$ is a Borel measurable function such
that $u(t,\cdot)\in \calg^1(\bfC,\R)$ for every $t\in[0,T]$
and
\begin{equation}\label{iposuu}
   |u(t,x)|+|\nabla_x u(t,x)|\leq C(1+|x|)^{m},
\end{equation}
for some $C>0, m\geq 0$ and
   for every $t\in [0,T]$, $x\in E$.

Then for every $x\in \bfC$, $i=1,\ldots,d$ and $0\le t\leq T'<T$   the processes
$\{u(s,X_s^{t,x})$, $s\in [t,T]\}$ and
 $W^i$ admit a joint quadratic variation
on the interval  $[t,T']$, given by the formula:
$$
\<u(\cdot,X_\cdot^{t,x}),W^i\>_{[t,T']} =\int_t^{T'}
\sigma^i(s,X_s^{t,x})\cdot\nabla_0 u(s,X_s^{t,x})\;ds
=\sum_{k=1}^n\int_t^{T'} \sigma^i_k(s,X_s^{t,x})\cdot\nabla_0^k
u(s,X_s^{t,x})\;ds.
$$
\end{theorem}

\noindent {\bf Proof.} For the sake of simplicity we write the
proof in the case $t=0$, the general case being deduced by the
same arguments.

We  fix $x\in \bfC$, $T'\in (0,T)$ and we denote $X^{0,x}$ by $X$ for
simplicity. Thus, $X_t=y(t+\cdot)$, $t\in [0,T]$, satisfies
$$
dy(t)=b(t,X_t)\;dt + \sigma(t,X_t)\;dW_t, \qquad X_0=x.
$$
We will use the results on the Malliavin derivatives stated in
Theorem \ref{derdimalliavinSDDE}, and in particular formula (\ref{eqdermallavinSDDE})
that, in view of (\ref{formulachainruleperlasol}), can be written in the form:
\begin{equation}\label{formuladerivmall}
    D_sy(t)=\sigma(s,X_s)+ \int_s^tD_s[b(r,X_r)]\;dr
+\int_s^tD_s[\sigma(r,X_r)]\;dW_r,
\end{equation}
for $0\le t\le s\le T$. Noting that $\nabla_x b(t,x)$ and $\nabla_x\sigma(t,x)$
are bounded by the
Lipschitz constant $L$ of $b(t,\cdot)$ and $\sigma(t,\cdot)$, it follows from
(\ref{formulachainruleperlasol}) that for every $r\in [0,T]$
\begin{equation}\label{stimedermall}
    \|D_\cdot [b(r,X_r)]\|^2\le L^2\int_0^T\sup_{t\in [s,T]}|D_s y(t)|^2\,ds,
\qquad \|D_\cdot [\sigma(r,X_r)]\|^2\le L^2\int_0^T\sup_{t\in [s,T]}|D_s y(t)|^2\,ds,
\end{equation}
where $\|\cdot\|$ denotes the norm in $L^2([0,T];\R^d)$.

We have to prove that
$$
\begin{array}{l}\dis
C^\epsilon:=
C^\epsilon_{[0,T']}(u(\cdot,X_\cdot),W^i)=\frac{1}{\epsilon}
\int_0^{T'} (u(t+\epsilon,X_{t+\epsilon})- u(t,X_{t}))
(W^i_{t+\epsilon}-W^i_t)\;dt
\\\dis\qquad\qquad
\to \int_0^{T'} \sigma^i(t,X_t)\cdot \nabla_0 u(t,X_t)\;dt
\end{array}
$$
in probability, as $\epsilon\to 0$.

We need to re-write $C_\epsilon$ in an appropriate way, fixing $
\epsilon>0$ so small that $T'+ \epsilon\leq T$. We first explain
our argument by writing down some informal passages: by the rules
of Malliavin calculus  we have, for a.a. $t \in [0,T']$,
\begin{equation}\label{primadagiustif}
  \begin{array}{l}\dis
(u(t +\epsilon,X_{t +\epsilon})- u(t ,X_{t })) (W^i _{t
+\epsilon}-W^i _t )= (u(t +\epsilon,X_{t +\epsilon})- u(t ,X_{t
}))e_i^* \int^{t +\epsilon}_{t }dW_s
\\\dis\quad
= \int^{t +\epsilon}_{t } D_s^i(u(t +\epsilon,X_{t +\epsilon})-
u(t ,X_{t }))\; ds +\int^{t +\epsilon}_{t } (u(t +\epsilon,X_{t
+\epsilon})- u(t ,X_{t }))e_i^* \hat{d}W_s,
\end{array}
\end{equation}
where the symbol $\hat{d}W$ denotes the Skorohod integral, and by
$e_i$ we denote the $i$-th component of the canonical basis of
$\R^d$ and by $e_i^*$ its transpose (row) vector. Integrating over
$[0,T']$ with respect to $t $ and interchanging integrals gives
\begin{equation}\label{secondadagiustif}
  \begin{array}{lll}\dis
\epsilon \;C^\epsilon &=&\dis \int_0^{T'}\int^{t +\epsilon}_{t }
D_s^i(u(t +\epsilon,X_{t +\epsilon})- u(t ,X_{t }))\; ds\;dt
\\&&\dis
+\int_0^{T'+\epsilon} \int_{(s-\epsilon)^+}^{s\wedge T'} (u(t
+\epsilon,X_{t +\epsilon})- u(t ,X_{t }))\;dt \; e_i^*\;
\hat{d}W_s.
\end{array}
\end{equation}
To justify (\ref{primadagiustif}) and (\ref{secondadagiustif})
rigorously we proceed as follows. To shorten notation we define
$$
v_t =(u(t +\epsilon,X_{t +\epsilon})- u(t ,X_{t }))\; 1_{[0,T']}(t
),\qquad t \in [0,T],
$$
$$
A^\epsilon=\{(t ,s)\in [0,T]\times [0,T]\;:\; 0\leq t \leq T', t
\leq s\leq t +\epsilon\}.
$$
Using Corollary \ref{chainruleperlasol} and formula (\ref{stimalpderivmall}), it is easy to show
 that, for all $t $, $v_t $ belongs to $\D^{1,2}$ and the
process $v_t \, 1_{A^\epsilon}(t ,\cdot)$ belongs to
$L^2(\Omega\times [0,T])$. By
 \cite{NuPa} Theorem 3.2 (see also
\cite{Nu} Section 1.3.1 (4)) we conclude that $v_t
\,1_{A^\epsilon}(t ,\cdot)\,e_i^*$ is Skorohod integrable and the
formula
\begin{equation}\label{eintegrsecondosko}
  \int_0^Tv_t \, 1_{A^\epsilon}(t ,s)\,e_i^*\;  \hat{d}W_s=
v_t\, \int_0^T 1_{A^\epsilon}(t ,s)\;e_i^*\;  \hat{d}W_s
-\int_0^TD_s^iv_t \, 1_{A^\epsilon}(t ,s)\; ds =:z_t ,
\end{equation}
holds provided $z_t $ belongs to $L^2(\Omega)$.
Since $\int_0^T 1_{A^\epsilon}(t ,s)\;  \hat{d}W_s$ coincides with
the Ito integral $\int_0^T 1_{A^\epsilon}(t ,s)\;  dW_s = (W_{t
+\epsilon}-W_t ) 1_{[t,T']}(t )$, it is in fact easy to verify
that we even have $z\in L^2(\Omega\times [0,T])$; thus
(\ref{eintegrsecondosko}) holds  for a.a. $t $, and
(\ref{eintegrsecondosko}) yields (\ref{primadagiustif}) for a.a.
$t\in [0,T']$.

Next we wish to show that the process $ \int_0^Tv_t
1_{A^\epsilon}(t,\cdot)\; dt\,e_i$ is Skorohod integrable and to
compute its integral, which occurs in the right-hand side of
(\ref{secondadagiustif}). For  arbitrary $G\in \D^{1,2}$, by the
definition of the Skorohod integral and by
(\ref{eintegrsecondosko}),
$$
  \begin{array}{lll}\dis
\E\int_0^T\left\<\int_0^Tv_t  1_{A^\epsilon}(t ,s) \; dt\,e_i
,D_sG\right\>_{\R^d}\;ds &=&\dis \int_0^T\E\int_0^T\left\<v_t
1_{A^\epsilon}(t ,s)\,e_i, D_sG\right\>_{\R^d}\;ds\; dt
\\&=&\dis
\int_0^T\E\left[G\int_0^Tv_t  1_{A^\epsilon}(t ,s) \;e_i^*\;
\hat{d}W_s\right]\; dt
\\&=&\dis
\E\left[G\int_0^Tz_t \; dt \right].
\end{array}
$$
This shows, by definition, that $ \int_0^Tv_t  1_{A^\epsilon}(t
,\cdot)\; dt\,e_i $ is Skorohod integrable and
$$
\int_0^T \int_0^Tv_t  1_{A^\epsilon}(t ,s)\; dt\;e_i^* \;
\hat{d}W_s = \int_0^Tz_t \; dt = \int_0^T
  \int_0^Tv_t  1_{A^\epsilon}(t ,s)\;  e_i^*\;\hat{d}W_s
\; dt .
$$
 Recalling
(\ref{eintegrsecondosko}) we obtain
$$
\int_0^T \int_0^Tv_t  1_{A^\epsilon}(t ,s)\; dt \;
e_i^*\;\hat{d}W_s = \int_0^Tv_t (W^i_{t +\epsilon}-W^i_t )\;
1_{[t,T']}(t )\; dt -\int_0^T \int_0^TD_s^iv_t \,e_i\,
1_{A^\epsilon}(t ,s)\; ds\; dt ,
$$
and (\ref{secondadagiustif}) is proved.

Recalling that $D_s(u(t ,X_{t }))=0$
for $s>t $ by adaptedness, and using the chain rule (\ref{formulachainruleperlasol}) for the
Malliavin derivative we have, for a.a. $s,t $ with $s\in [t ,t
+\epsilon]$,
$$
D_s(u(t +\epsilon,X_{t +\epsilon})- u(t ,X_{t }))=
D_s(u(t +\epsilon,X_{t +\epsilon})) = \int_{[-r,0]} D_sy(t
+\epsilon+\theta)\cdot \nabla_x u(t +\epsilon,X_{t
+\epsilon})(d\theta)
$$
 and from (\ref{secondadagiustif}) we deduce
$$\begin{array}{lll}\dis
C^\epsilon &=&\dis \frac{1}{\epsilon} \int_0^{T'}\int^{t
+\epsilon}_{t } \int_{[-r,0]} D_s^iy(t +\epsilon+\theta)\cdot
\nabla_x u(t +\epsilon,X_{t +\epsilon})(d\theta)
\; ds\;dt\\
&+&\dis \frac{1}{\epsilon} \int_0^{T'+\epsilon}
\int_{(s-\epsilon)^+}^{s\wedge T'} (u(t +\epsilon,X_{t
+\epsilon})- u(t ,X_{t }))\;dt \; e_i^*\; \hat{d}W_s
\\
&=:&\dis I_1^\epsilon+I_2^\epsilon.
\end{array}
$$

Now we let $\epsilon\to 0$, and we first claim that
$I_2^\epsilon\to 0$ in probability. To prove this, it is enough to show that the
process $\frac{1}{\epsilon} \int_0^T(u(t +\epsilon,X_{t
+\epsilon})- u(t ,X_{t }))\;1_{A^\epsilon}(t ,\cdot)\; dt $
converges to $0$ in $\L^{1,2}$. Indeed, since the Skorohod
integral is a bounded linear operator from $\L^{1,2}$ to
$L^2(\Omega)$, this  implies that
$$
I_2^\epsilon= \int_0^T \frac{1}{\epsilon} \int_0^T(u(t
+\epsilon,X_{t +\epsilon})- u(t ,X_{t })) 1_{A^\epsilon}(t ,s)\;
dt \; e_i^*\;\hat{d}W_s \to 0
$$
in $L^2(\Omega)$.
 We prove, more generally, that for an
arbitrary element $y\in \L^{1,2}(\R)$, if we set
$$
T^\epsilon(y)_s= \frac{1}{\epsilon} \int_0^T(y_{t +\epsilon}- y_{t
})\;1_{A^\epsilon}(t ,s)\; dt = \frac{1}{\epsilon}
\int_{(s-\epsilon)\vee t}^{s\wedge T'} (y_{t +\epsilon}- y_{t
})\;dt ,\qquad s\in [0,T],
$$
then the process $T^\epsilon(y)$ converges to $0$ in
$\L^{1,2}(\R)$. Let us recall that $\L^{1,2}(\R)$ is isomorphic to
$L^2([0,T]; \D^{1,2}(\R))$. It is clear that $T^\epsilon(y)\to 0$ in $\L^{1,2}(\R)$
if $y$ belongs to $C([0,T]; \D^{1,2}(\R))$, a dense subspace of
$L^2([0,T]; \D^{1,2}(\R))$. So to prove the claim it is enough to
show that the norm of $T^\epsilon$, as an operator on
$\L^{1,2}(\R)$, is bounded uniformly with respect to $\epsilon$.
We have
$$
\begin{array}{lll}\dis
|T^\epsilon(y)_s|^2_{\D^{1,2}(\R)} &\le&\dis \frac{1}{\epsilon^2}
\int_0^T1_{A^\epsilon}(t ,s)\; dt \int_0^T|y_{t +\epsilon}- y_{t
}|^2_{\D^{1,2}(\R)}\;1_{A^\epsilon}(t ,s)\; dt
\\&\le&\dis
\frac{1}{\epsilon} \int_0^T|y_{t +\epsilon}- y_{t
}|^2_{\D^{1,2}(\R)}\;1_{A^\epsilon}(t ,s)\; dt ,
\end{array}
$$
$$
\begin{array}{lll}\dis
|T^\epsilon(y)|^2_{\L^{1,2}(\R)} &=&\dis
\int_0^T|T^\epsilon(y)_s|^2_{\D^{1,2}(\R)}\;ds
\\
&\le&\dis \frac{1}{\epsilon} \int_0^T|y_{t +\epsilon}- y_{t
}|^2_{\D^{1,2}(\R)} \int_0^T1_{A^\epsilon}(t ,s)\;ds\; dt
\\
\\
&\le&\dis \int_0^{T'}|y_{t +\epsilon}- y_{t }|^2_{\D^{1,2}(\R)}\;
dt
\\
&\le&\dis 2|y|^2_{\L^{1,2}(\R)}.
\end{array}
$$
This shows the required bound, and completes the proof that
$I_2^\epsilon\to 0$ as $\epsilon\to 0$.

Now we proceed to compute the limit of $I_1^\epsilon$.
 We note that, by adaptedness,
 $D_sy(t +\epsilon+\theta)=0$ for $s>t +\epsilon+\theta$,
 so that
$$
I_1^\epsilon= \frac{1}{\epsilon} \int_0^{T'}\int^{t +\epsilon}_{t
} \int_{[s-t-\epsilon,0]} D_s^iy(t +\epsilon+\theta)\cdot \nabla_x
u(t +\epsilon,X_{t +\epsilon})(d\theta) \; ds\;dt.
$$
For fixed $t$, let us exchange integrals with respect to $ds$ and
$\nabla_x u(t +\epsilon,X_{t +\epsilon})(d\theta)$ obtaining
$$
I_1^\epsilon= \frac{1}{\epsilon} \int_0^{T'} \int_{[-\epsilon,0]}
\int^{t +\epsilon+\theta}_{t } D_s^iy(t +\epsilon+\theta)\;
ds\cdot \nabla_x u(t +\epsilon,X_{t +\epsilon})(d\theta) \;dt.
$$
Next we replace   $D_sy(t +\epsilon+\theta)$ by the expression
given by (\ref{formuladerivmall}) and we obtain
$$\begin{array}{l}
I_1^\epsilon=\dis \frac{1}{\epsilon} \int_0^{T'}
\int_{[-\epsilon,0]} \int^{t +\epsilon+\theta}_{t }
\sigma^i(s,X_s)\; ds\cdot \nabla_x u(t +\epsilon,X_{t
+\epsilon})(d\theta) \;dt
\\\dis\qquad
+ \frac{1}{\epsilon} \int_0^{T'} \int_{[-\epsilon,0]} \int^{t
+\epsilon+\theta}_{t } \int_s^{t
+\epsilon+\theta}D^i_s[b(r,X_r)]\;dr \; ds\cdot \nabla_x u(t
+\epsilon,X_{t +\epsilon})(d\theta) \;dt
\\\dis\qquad
+ \frac{1}{\epsilon} \int_0^{T'} \int_{[-\epsilon,0]} \int^{t
+\epsilon+\theta}_{t } \int_s^{t
+\epsilon+\theta}D^i_s[\sigma(r,X_r)]\;dW_r \; ds\cdot \nabla_x
u(t +\epsilon,X_{t +\epsilon})(d\theta) \;dt
\\\qquad
=: J_1^\epsilon+J_2^\epsilon+J_3^\epsilon.
\end{array}
$$

We first show that $J_3^\epsilon\to 0$ in $L^1(\Omega)$.  Since,
by (\ref{iposuu}),
$$
|\nabla_x u(t +\epsilon,X_{t +\epsilon})|\le
C(1+\sup_{t\in [0,T]}|X_t|_\bfC)^m,
$$
then, using the notation
$|\nabla_x u(t +\epsilon,X_{t +\epsilon})|(d\theta)$ to indicate
the total variation measure, we have
$$
\begin{array}{l}\dis
|J_3^\epsilon|\le \frac{1}{\epsilon} \int_0^{T'}
\int_{[-\epsilon,0]} \int^{t +\epsilon+\theta}_{t } \left|
\int_s^{t +\epsilon+\theta}D_s^i[\sigma(r,X_r)]\;dW_r \right| \;
ds\; |\nabla_x u(t +\epsilon,X_{t +\epsilon})|(d\theta) \;dt,
\\\dis
\qquad \le C(1+\sup_{t\in [0,T]}|X_t|_\bfC)^m \frac{1}{\epsilon}
\int_0^{T'} \sup_{\theta\in[-\epsilon,0]} \int^{t
+\epsilon+\theta}_{t } \left| \int_s^{t
+\epsilon+\theta}D_s^i[\sigma(r,X_r)]\;dW_r \right| \; ds\; \;dt,
\\\dis
\qquad \le C(1+\sup_{t\in [0,T]}|X_t|_\bfC)^m \frac{1}{\epsilon}
\int_0^{T'} \int^{t +\epsilon}_{t }
\sup_{\theta\in[-\epsilon,0]}\left| \int_s^{t
+\epsilon+\theta}D_s^i[\sigma(r,X_r)]\;dW_r \right| \; ds\; \;dt.
\end{array}
$$
Taking the $L^1(\Omega)$ norm of both sides and using the H\"older and the Doob
maximal inequality we have
$$\begin{array}{l}
\|J_3^\epsilon\|_{L^1(\Omega)}
\\\dis
\quad \le C\|(1+\sup_{t\in [0,T]}|X_t|_\bfC)^m\|_{L^2(\Omega)}
\frac{1}{\epsilon} \int_0^{T'} \int^{t +\epsilon}_{t } \left\|
\sup_{\theta\in[-\epsilon,0]}\left| \int_s^{t
+\epsilon+\theta}D_s^i[\sigma(r,X_r)]\;dW_r
\right|\right\|_{L^2(\Omega)} \; ds\; \;dt
\\\dis
\quad \le \frac{C}{\epsilon} \int_0^{T'} \int^{t +\epsilon}_{t }
\left\| \int_s^{t +\epsilon}D_s^i[\sigma(r,X_r)]\;dW_r
\right\|_{L^2(\Omega)} \; ds\; \;dt
\\\quad=\dis
\frac{C}{\epsilon} \int_0^{T'} \int^{t +\epsilon}_{t } \left(
\int_s^{t +\epsilon}\E\, \left|D_s^i[\sigma(r,X_r)]
\right|^2\;dr\right)^{1/2} \; ds\; \;dt.
\end{array}
$$
Denoting for simplicity $h(s,r)=\E\, \left|D_s^i[\sigma(r,X_r)]
\right|^2$ we obtain
$$\begin{array}{lll}\dis
\|J_3^\epsilon\|_{L^1(\Omega)} &\le& \dis \frac{C}{\epsilon}
\int_0^{T'} \int^{t +\epsilon}_{t } \left(\int^{t +\epsilon}_{t}
h(s,r)\;dr \right)^{1/2} \; ds\;dt
\\&=&\dis
\frac{C}{\sqrt{\epsilon}} \int_0^{T'} \left(\int^{t +\epsilon}_{t
} \int^{t +\epsilon}_{t} h(s,r)\;dr\; ds \right)^{1/2} \;dt
\\&\le&\dis
{C}{\sqrt{T'}} \left(\int_0^{T'}\left[ \frac{1}{\epsilon} \int^{t
+\epsilon}_{t } \int^{t +\epsilon}_{t} h(s,r)\;dr\; ds\right]\;dt
\right)^{1/2}.
\end{array}
$$
Let us note that $h\in L^1([0,T]^2)$, since by
(\ref{stimedermall}) we have
$$
 \int^{T}_{0}\int^{T}_{0} h(s,r)\;dr \; ds \le\E\,\int_0^T
\|D_\cdot [\sigma(r,X_r)]\|^2 dr
\le
L^2T\int_0^T\E\,\sup_{t\in [s,T]}|D_s y(t)|^2ds <\infty.
$$
Let us define the operator $\cala_\epsilon:L^1([0,T]^2)\to
L^1([0,T])$ by
$$
(\cala_\epsilon k)(t)= \frac{1}{\epsilon} \int^{(t
+\epsilon)\wedge T}_{t } \int^{(t +\epsilon)\wedge T}_{t}
k(s,r)\;dr\; ds, \qquad k\in L^1([0,T]^2).
$$
Then we have $\|J_3^\epsilon\|_{L^1(\Omega)} \le {C}{\sqrt{T'}}
\|\cala_\epsilon h\|_{L^1([0,T])}^{1/2}$, so to prove that
$J_3^\epsilon\to 0$ in ${L^1(\Omega)}$ it is enough to show that
$\cala_\epsilon k\to 0$ in $L^1([0,T])$ for every $k \in
L^1([0,T]^2)$. This is obvious if $k$ is in the space of bounded
functions on $[0,T]^2$, a dense subspace of $L^1([0,T]^2)$. So it
is enough to show that $ \|\cala_\epsilon k\|_{L^1([0,T])}\le C
\|k\|_{L^1([0,T]^2)}$ for some constant $C$ and for every $k \in
L^1([0,T]^2)$. This follows from the inequalities
$$
\begin{array}{l}\dis
\int_0^T|\cala_\epsilon k(t)|\;dt \le
\frac{1}{\epsilon}\int_0^T\int_0^T\int_0^T
|k(s,r)|1_{t<s<(t+\epsilon)\wedge T}1_{t<r<(t+\epsilon)\wedge T}
dr\;ds\;dt
\\\dis\qquad
= \frac{1}{\epsilon}\int_0^T\int_0^T|k(s,r)|\left[\int_0^T
1_{(s-\epsilon)^+<t< s}1_{(r-\epsilon)^+<t< s} dt\right]\;ds\;dr
\le \int_0^T\int_0^T|k(s,r)|\;ds\;dr,
\end{array}
$$
since the term in square brackets is less or equal to $\epsilon$.

This finishes the proof that $J_3^\epsilon\to 0$ in $L^1(\Omega)$,
hence in probability. In a similar and simpler way one proves that
$J_2^\epsilon\to 0$ in probability.

To finish the proof of the proposition it remains to compute the
limit of $J_1^\epsilon$. Exchanging integrals with respect to $ds$
and $\nabla_x u(t +\epsilon,X_{t +\epsilon})(d\theta)$, and then
using another change of variable we have
$$\begin{array}{l}
J_1^\epsilon=\dis \frac{1}{\epsilon} \int_0^{T'}
\int_{[-\epsilon,0]} \int^{t +\epsilon+\theta}_{t }
\sigma^i(s,X_s)\; ds\cdot \nabla_x u(t +\epsilon,X_{t
+\epsilon})(d\theta) \;dt
\\\dis\qquad
=\frac{1}{\epsilon} \int_0^{T'} \int^{t +\epsilon}_{t }
\int_{[s-t-\epsilon,0]} \sigma^i(s,X_s)\cdot \nabla_x u(t
+\epsilon,X_{t +\epsilon})(d\theta) \;ds\;dt
\\\dis\qquad
=\frac{1}{\epsilon} \int_0^{T'} \int^{t +\epsilon}_{t }
\sigma^i(s,X_s)\cdot \nabla_x u(t +\epsilon,X_{t
+\epsilon})([s-t-\epsilon,0]) \;ds\;dt
\\\dis\qquad
=\frac{1}{\epsilon} \int_\epsilon^{T'+\epsilon} \int_{t
-\epsilon}^{t} \sigma^i(s,X_s)\cdot \nabla_x u(t,X_{t})([s-t,0])\;
ds \;dt
\\\dis\qquad
=\frac{1}{\epsilon} \int_\epsilon^{T'+\epsilon} \int_{t
-\epsilon}^{t} \sigma^i(t,X_t)\cdot \nabla_x u(t,X_{t})([s-t,0])\;
ds \;dt
\\\dis\qquad
+\frac{1}{\epsilon} \int_\epsilon^{T'+\epsilon} \int_{t
-\epsilon}^{t} \{\sigma^i(s,X_s)-\sigma^i(t,X_t)\}\cdot \nabla_x
u(t,X_{t})([s-t,0])\; ds \;dt
\\\dis\qquad
=: H_1^\epsilon+H_2^\epsilon.
\end{array}
$$

Next we show that $H_2^\epsilon\to 0$, $\P$-a.s. Since
$$
|\nabla_x u(t +\epsilon,X_{t +\epsilon})|\le
C(1+\sup_{t\in [0,T]}|X_t|_\bfC)^m,
$$
we have
$$\begin{array}{lll}
|H_2^\epsilon|&\le&\dis C(1+\sup_{t\in [0,T]}|X_t|_\bfC)^m
\frac{1}{\epsilon} \int_\epsilon^{T'+\epsilon} \int_{t
-\epsilon}^{t} |\sigma^i(s,X_s)-\sigma^i(t,X_t)|\; ds \;dt
\\&\le &\dis
C(1+\sup_{t\in [0,T]}|X_t|_\bfC)^m \int_0^{T}\frac{1}{\epsilon}
\int_{(t -\epsilon)^+}^{t} |\sigma^i(s,X_s)-\sigma^i(t,X_t)|\; ds
\;dt.
\end{array}
$$
Let us fix $\omega\in\Omega$ and note that, $\P$-a.s., $\sigma^i(\cdot,X_\cdot)\in
L^1([0,T])$. Let us define the operator
$\calb_\epsilon:L^1([0,T])\to L^1([0,T])$ as
$$
(\calb_\epsilon k)(t)= \frac{1}{\epsilon} \int_{(t
-\epsilon)^+}^{t} |k(s)-k(t)|\; ds, \qquad k\in L^1([0,T]).
$$
Then we have $|H_2^\epsilon| \le C(1+\sup_{t\in [0,T]}|X_t|_\bfC)^m\|\calb_\epsilon
\sigma^i(\cdot,X_\cdot) \|_{L^1([0,T])}$, $\P$-a.s., so to prove
that $H_2^\epsilon\to 0$ in probability it is enough to show that
$\calb_\epsilon k\to 0$ in $L^1([0,T])$ for every $k \in
L^1([0,T])$. This is obvious if $k$ is in the space of continuous
functions on $[0,T]$, a dense subspace of $L^1([0,T])$. So it is
enough to show that $ \|\calb_\epsilon k\|_{L^1([0,T])}\le C
\|k\|_{L^1([0,T])}$ for some constant $C$ and for every $k \in
L^1([0,T])$. This follows from the inequality
$$
|(\calb_\epsilon k)(t)|\le \frac{1}{\epsilon} \int_{(t
-\epsilon)^+}^{t} |k(s)| \; ds +|k(t)|,
$$
which implies
$$\begin{array}{l}
\dis \int_0^T|(\calb_\epsilon k)(t)|\;dt\le \|k\|_{L^1([0,T])}+
\frac{1}{\epsilon}\int_0^T \int_{(t -\epsilon)^+}^{t} |k(s)| \; ds
\;dt
\\\dis\qquad
=\|k\|_{L^1([0,T])}+ \frac{1}{\epsilon}\int_0^T \int^{(s
+\epsilon)\wedge T}_{s} |k(s)| \;dt\; ds
\\\dis\qquad
\le \|k\|_{L^1([0,T])}+ \int_0^T |k(s)| \; ds =
2\|k\|_{L^1([0,T])}.
\end{array}
$$
This finishes the proof that $H_2^\epsilon\to 0$ $\P$-a.s., hence
in probability.

It remains to consider the term
$$\begin{array}{l}\dis
H_1^\epsilon=\frac{1}{\epsilon} \int_\epsilon^{T'+\epsilon}
\sigma^i(t,X_t)\cdot\int_{t -\epsilon}^{t} \nabla_x
u(t,X_{t})([s-t,0])\; ds \;dt
\\\dis\qquad
= \frac{1}{\epsilon} \int_\epsilon^{T'+\epsilon}
\sigma^i(t,X_t)\cdot\int_{t -\epsilon}^{t} \int_{[s-t,0]} \nabla_x
u(t,X_{t})(d\theta)\; ds \;dt
\\\dis\qquad
= \frac{1}{\epsilon} \int_\epsilon^{T'+\epsilon}
\sigma^i(t,X_t)\cdot \int_{[-\epsilon,0]} \int_{t
-\epsilon}^{t+\theta}\; ds \nabla_x u(t,X_{t})(d\theta) \;dt
\\\dis\qquad
= \frac{1}{\epsilon} \int_\epsilon^{T'+\epsilon}
\sigma^i(t,X_t)\cdot \int_{[-\epsilon,0]} (\theta+\epsilon)\;
\nabla_x u(t,X_{t})(d\theta) \;dt
\\\dis\qquad
= \int_\epsilon^{T'+\epsilon} \sigma^i(t,X_t)\cdot \int_{[-r,0]}
\left(1+\frac{\theta}{\epsilon}\right)^+\; \nabla_x
u(t,X_{t})(d\theta) \;dt.
\end{array}
$$
We clearly have, $\P$-a.s.,
$$
\begin{array}{lll}\dis
\int_{[-r,0]} \left(1+\frac{\theta}{\epsilon}\right)^+ \nabla_x
u(t,X_{t})(d\theta) & \to &\dis \int_{[-r,0]} 1_{\{0\}}(\theta)\;
\nabla_x u(t,X_{t})(d\theta)
\\
&=&\dis \nabla_x u(t,X_{t})(\{0\})= \nabla_0 u(t,X_{t}),
\end{array}
$$
and by dominated convergence,  $\P$-a.s.,
$$
H_1^\epsilon\to \int_0^{T'} \sigma^i(t,X_t)\cdot \nabla_0
u(t,X_{t}) \;dt.
$$

This shows that $C^\epsilon$ converges in probability and its
limit is
$$
\<u(\cdot,X_\cdot),W^i\>_{[0,T']} =\int_0^{T'}
\sigma^i(t,X_t)\cdot \nabla_0 u(t,X_t)\;dt. \qed
$$

\section{The forward-backward system with delay}
\label{sezioneforback}

In this section we will
 discuss existence, uniqueness and regular
dependence on the initial data of the following
forward-backward system: for given $t\in [0,T]$ and $x\in\bfC=C([-r,0];\R^n)$,
\begin{equation}\label{fbsde}
    \left\{\begin{array}{l}\dis dy_\tau =
b(\tau,{X}_\tau)\; d\tau
+\sigma(\tau,{X}_\tau)\; dW_\tau,\quad \tau\in
[t,T]\subset [0,T],
\\\dis
X_t=x,
\\\dis
 dY_\tau=\psi(\tau,X_\tau,Y_\tau,Z_\tau)\;d\tau+Z_\tau\;dW_\tau,
  \\\dis
  Y_T=\phi(X_T),
\end{array}\right.
\end{equation}
Here we use the notation $X_\tau(\theta)=y_{\tau+\theta}$, $ \theta\in
[-r,0]$, as before, so the first equation in (\ref{fbsde})
is the same as (\ref{eqscalare-s}). We extend the definition of $X$ setting
$X_s=x$ for $0\leq s\leq t$. The second equation in
(\ref{fbsde}), namely
\begin{equation}\label{bsde}
    \left\{\begin{array}{l}\dis
 dY_\tau=\psi(\tau,X_\tau,Y_\tau,Z_\tau)\;d\tau+Z_\tau\;dW_\tau,
 \qquad \tau\in [0,T],
  \\\dis
  Y_T=\phi(X_T),
\end{array}\right.
\end{equation}
is of backward type.
Under suitable assumptions on the coefficients
 $\psi:[0,T]\times\bfC\times\R\times\R^d
\rightarrow\mathbb{R}$ and $\mathbb{\phi}:\bfC\rightarrow\mathbb{R}$
we will look for a solution consisting of a pair  $  $ of predictable processes,
taking values in $\mathbb{R}\times\R^d$, such that $Y$ has
continuous paths and
\[
\|\left(  Y,Z\right)\|^2_{\mathbb{K}_{cont}}:=
\mathbb{E}\sup_{\tau\in\left[  0,T\right]  }\left\vert Y_{\tau}\right\vert
^{2}+\mathbb{E}\int_{0}^{T}\left\vert Z_{\tau}\right\vert ^{2}d\tau<\infty,
\]
see e.g. \cite{PaPe1}. In the following we denote by $\mathbb{K}_{cont}\left(
\left[  0,T\right]  \right)$ the space of such processes.

The solution of (\ref{fbsde}) will be denoted by $(X_\tau, Y_\tau, Z_\tau)_{\tau\in[0,T]}$, or,
to stress the dependence on the initial time $t$ and on the
initial datum $x$, by $(X_\tau^{t,x}, Y_\tau^{t,x}, Z_\tau^{t,x})_{\tau\in[0,T]}$.

\begin{hypothesis}
\label{ip su psi}The maps $\psi:[0,T]\times\bfC\times\R\times\R^d
\rightarrow\mathbb{R}$ and $\mathbb{\phi}:\bfC\rightarrow\mathbb{R}$ are Borel
measurable and satisfy the following assumptions:
\begin{enumerate}
\item there exists $L>0$ such that
\begin{align*}
&\left\vert \psi\left( t,x,y,z_{1}\right)  -\psi\left(t,x,y,z_{2}\right)  \right\vert
\leq L\left\vert z_{1}-z_{2}\right\vert,\\ \nonumber
&\left\vert \psi\left(t,  x,y_1,z\right)  -\psi\left(t,x,y_2,z\right)  \right\vert
\leq L\left\vert y_{1}-y_{2}\right\vert,\\ \nonumber
\end{align*}
for every $t\in [0,T]$, $x\in \bfC$, $y,y_1,y_2 \in\R$
and $z, z_{1},z_{2}\in\R^d$;

\item  $\psi(
  t ,\cdot, \cdot , \cdot )  \in \calg^{1}\left( \bfC\times\R\times\R^d,\R\right)  $
  for every $t\in\left[  0,T\right]  $;

\item there exist $K>0$ and $m\geq0$ such that
\[
\left\vert \nabla_{x}\psi\left(  t,x,y,z\right) \right\vert \leq
K\left(  1+\left\vert x\right\vert _{\bfC}+\vert y\vert\right)
^{m}\left(  1+\left\vert z\right\vert\right)
\]
for every $t\in\left[  0,T\right]  $, $x\in \bfC$, $y\in\R$ and $z\in\R^d$;

\item $\phi\in \calg^{1}\left(  \bfC,\R\right)  $ and there exist $K>0$ and $m\geq0$ such that
\[
\left\vert \nabla\phi(x)  \right\vert \leq
K\left(  1+\left\vert x\right\vert _{\bfC}\right)^{m},\qquad x\in \bfC.
\]
\end{enumerate}
\end{hypothesis}

Under these assumptions we can state a result on existence and uniqueness of a
solution of the forward-backward system (\ref{fbsde})
and on its regular dependence on $x$.

\begin{proposition}\label{propBSDE}
Assume that Hypotheses \ref{ipotesi1} and \ref{ip su psi} hold
true. Then the forward-backward system (\ref{fbsde}) admits a
unique solution $\left(X^{t,x},Y^{t,x},Z^{t,x}\right)  \in
\cals^p([0,T];\bfC)\times \mathbb{K}_{cont}\left( \left[
0,T\right]  \right)  $ for every $\left(  t,x\right)\in [0,T]\times \bfC$. Moreover,
the map $\left(  t,x\right)
\mapsto( X^{t,x}, Y^{t,x},Z^{t,x})  $ belongs to the space
 $\calg^{1}\left([0,T]\times \bfC,\cals^p([0,T];\bfC)\times\mathbb{K}_{cont}\left( \left[
0,T\right]  \right)  \right) $. Finally, the following estimate
holds
true: for every $p\geq2$ there exists $C>0$ such that
\[
\left[ \mathbb{E}\sup_{\tau\in\left[  0,T\right]  }\left\vert \nabla_{x}Y_\tau^{t,x}
  \right\vert ^{p}\right]  ^{1/p}\leq C
\left(  1+\left\vert x\right\vert _{\bfC}^{\left(  m+1\right)  ^{2}}\right),
\qquad t\in [0,T],x\in \bfC.
\]
\end{proposition}

\noindent {\bf Proof.}
We only give a sketch of the proof. The forward equation has a unique
solution by Theorem \ref{teoSDDE}. Existence and uniqueness of the solution of the backward
equation follows  from the classical result \cite{PaPe1}.

In Theorem \ref{teoSDDE} we have shown that
the map $x\mapsto X^{t,x} $ belongs to
$ C^{1}\left( \bfC,\mathcal{S}^{p}\left(  \left[ 0,T\right]
;\bfC\right)  \right)  $  for every $2\leq p<\infty$.
Then  the proof of continuity and differentiability of
$\left(  t,x\right)
\mapsto( X^{t,x}, Y^{t,x},Z^{t,x})  $  in the appropriate norms, as well as the final
estimate on $\nabla_{x}Y_\tau^{t,x}$,
can be achieved as in  Proposition 5.2 in \cite{fute}.
 The only difference is that
in \cite{fute} the process $X^{t,x}$ takes values in a
Hilbert space, while in our context it takes values
in the Banach space $\bfC$; nevertheless the same arguments apply
 (see also
\cite{Mas} for a similar result in Banach spaces).
 \qed

\begin{corollary}\label{proprietadiv}
 Assume that Hypotheses \ref{ipotesi1} and \ref{ip su psi} hold
true. Then the function $v:[0,T]\times \bfC\to \R$
defined by
\begin{equation}\label{defdiv}
v(t,x)=Y_t^{t,x},\qquad
t\in[0,T],x\in \bfC,
\end{equation}
 belongs to
$\calg^{0,1}\left([0,T]\times  \bf C;\mathbb{R}
\right)  $. Moreover there exists $C>0$ such that
$$
\left\vert \nabla_{x}v\left(
t,x\right)  \right\vert \leq C\left(
1+\left\vert x\right\vert _{\bf C}^{\left(  m+1\right)  ^{2}}\right),
\qquad t\in\left[  0,T\right]  ,x\in \bf C.
$$
Finally, for every $ t\in\left[  0,T\right] $ and $x\in \bf C$,
 we have, $\P$-a.s,
\begin{equation}
 \label{identificazioneY}
Y^{t,x}_{s}=
v\left(  s,X^{t,x}_{s}\right), \qquad \text{\emph{for every }
} s\in [t,T],
\end{equation}
\begin{equation}
 \label{identificazioneZ}
Z^{t,x}_{s}=\nabla_0
v\left(  s,X^{t,x}_{s}\right) \sigma(s,X^{t,x}_{s}) , \qquad \text{\emph{for a.e.}
} s\in [t,T].
\end{equation}
\end{corollary}

\noindent {\bf Proof.} It is well known that $v(t,x)$ is deterministic, and
its properties are therefore a direct consequence of Proposition
\ref{propBSDE}.
Equality (\ref{identificazioneY}) is also a standard consequence of uniqueness
of the solution of the backward equation.

To prove  (\ref{identificazioneZ}) we consider the joint quadratic variation
of $Y^{t,x}$ and  the Wiener process
$W^i$
on
an interval $[t,T']$, with $T'<T$.
Taking into account the backward equation we obtain
$$
\< Y^{t,x},W^i\>_{[t,T']}=\int_t^{T'}Z^i_s\;ds.
$$
By Theorem \ref{esistejointvar} we have
$$
\<  v\left(  \cdot,X^{t,x}_{\cdot}\right),W^i\>_{[t,T']}=\int_t^{T'}
\sigma^i(s,X_s^{t,x})\cdot \nabla_0 v(s,X_s^{t,x})
\;ds,
$$
so that (\ref{identificazioneY}) implies  (\ref{identificazioneZ}).
\qed

\begin{remark}\begin{em}
If we strengthen slightly the regularity assumptions
and we require that, for all $t\in[0,T]$, the functions $b(t,\cdot),
\sigma(t,\cdot),\phi$ are continuously Fr\'echet differentiable
on $\bfC$ and $\psi(
  t ,\cdot, \cdot , \cdot )$ is  continuously Fr\'echet differentiable
on $ \bfC\times\R\times\R^d$, then
we can prove, with only minor changes in the proofs, that the function
$v$ defined in  (\ref{defdiv}) is  Fr\'echet differentiable
with respect to $x$ and the
 Fr\'echet derivative is a continuous function from
$[0,T]\times \bfC$  to the dual space $\bfC^*$
with respect to the usual norm (i.e. the variation norm).
\end{em}\end{remark}

\begin{remark}\begin{em}
\label{remark v} In the context of Proposition \ref{propBSDE}, the
law of the solution $(X^{t,x},Y^{t,x},Z^{t,x})$ is uniquely determined
by $,x$ and the coefficients $b,\sigma,\psi,\phi$. Since
 $v(t,x)$ is deterministic, hence determined by its law, we conclude that the function
 $v$
is a functional of the coefficients $b,\sigma,\psi,\phi$ and
does not depend on the particular choice of the probability space
$({\Omega},{\calf},{\mathbb{P}})  $
nor on the Wiener process
$  {W} $.
\end{em}\end{remark}

\section{Application to stochastic optimal control}
\label{sez-control}
\subsection{Strong formulation of the optimal control problem}

Let
$(  \Omega,\mathcal{F},\left(  \mathcal{F}_{t}\right)  _{t\geq0},\mathbb{P})  $
be a filtered probability space, satisfying the usual conditions, and let $W$ be
 an $\R^d$-valued standard Wiener process with respect to
 $(\mathcal{F}_{t})$ and $\mathbb{P}$.
We consider the following
controlled functional stochastic equation
on an interval $
[{t},T]\subset [0,T]$:
\begin{equation}
\label{SDDEcontrollata}
\left\{
\begin{array}{l}
dy^u_s =
b(s,y^u_{s+\cdot})\; ds
+\sigma(s,y^u_{s+\cdot})\,
[h(s,y^u_{s+\cdot},u_s)\;ds+\; dW_s]
,
\\
y^u_{t+\theta}=x(\theta),\quad \theta\in [-r,0],
\end{array}
\right.
\end{equation}
The coefficients $b$ and $\sigma$ satisfy the previous assumptions.
 $u(\cdot)$ denotes the
control and $y^u$ the corresponding solution. We assume that controls are
 $\left(  \mathcal{F}_{t}\right)  $-predictable process with values
in a given measurable space $(U,\calu)$.
The function
 $h:\left[
0,T\right]  \times \bfC \times U\to \R^{d}  $
is measurable and bounded.
We introduce again the
process
\begin{equation}\label{segmentou}
X^u_s=y^u_{s+\cdot}=\{y^u_{s+\theta},\; \theta\in [-r,0]\},
\qquad s\in [t,T],
\end{equation}
which now depends on the control and takes
values in  $\bfC = C([-r,0];\R^n)$,
so that equation (\ref{SDDEcontrollata}) can be rewritten as
\begin{equation}\label{SDDEcontrollata1}
\left\{
\begin{array}{l}
dy^u_s =
b(s,{X}^u_s)\; ds
+\sigma(s,X^u_s)\,
[h(s,X^u_s,u_s))\;ds+\; dW_s]
,\quad s\in
[{t},T],
\\
{X}_{{t}}={x}.
\end{array}
\right.
\end{equation}

We introduce the cost functional to minimize:
\begin{equation}\label{funzcostofinale}
J(t,x,u(\cdot))=
\E\int_t^Tg( u_s)\; ds
+\E\,
\phi(y^u_{T+\cdot})
=\E\int_t^Tg( u_s)\; ds
+
\E\,\phi(X^u_T),
\end{equation}
where $g:U\to [0,\infty)$ and $\phi:\bfC \rightarrow\R$  are given functions.

\begin{remark}\begin{em} Without any substantial change,
we could consider more general cost functionals of the form
\begin{equation}
\label{funzcosto1}
J(t,x,u(\cdot))=
\E\int_t^T\left[\ell(y^u_{s})+g(u_s)\right]\;ds+
\E\phi(y^u_{T+\cdot}),
\end{equation}
where $\ell:\R^n\rightarrow\R$. In fact, this kind of cost
can be put in the form (\ref{funzcostofinale}) as follows:
first note that
in equation (\ref{SDDEcontrollata})
we can assume  $r\geq T$, possibly extending the functions $b$
and $\sigma$ in the obvious way; next we define, for $x\in\bfC$,
\begin{equation*}
 \phi_0(x)=
\int_{t-T}^0\ell(x(s))\;ds.
\end{equation*}
so that
$ \phi_0(X^u_T)=\int_t^T\ell_0(y^u_s)\;ds$
and we conclude that
$$
J(t,x,u(\cdot))=
\E\int_t^Tg(u_s)\;ds
+\E [(\phi_0+
\phi)(X^u_T)],
$$
which has the required form.
In a similar way, under suitable assumptions,
one could consider even more general costs
of the form
$$
J(t,x,u(\cdot))=
\E\int_t^T\ell(s,y^u_{s},u_s)\;ds+
\E\,\phi(y^u_{T+\cdot}).
$$
However, we limit ourselves to
 cost functionals with the structure of
 (\ref{funzcostofinale}).
\end{em}
\end{remark}

To proceed further we need to introduce
 the hamiltonian function
$\psi:[0,T]\times \bfC\times\R^d\to\mathbb{R}$ defined,
for $t\in [0,T]$, $x\in \bfC,$ $z\in\R^d$, by the formula
\begin{equation}
\psi\left(t, x,z\right)  =\inf\left\{g(u)
+zh\left(t,x,u\right)  :u\in U\right\}   \label{hamiltoniana}
\end{equation}
and the corresponding, possibly empty, set of minimizers
\begin{equation}
\Gamma\left(t,x,z\right)  =\left\{u\in U,\,\;\, g( u)
+zh\left(t,x,u\right)   =\psi\left(t,x,z\right)  \right\}  .
\label{gamma}
\end{equation}

\begin{remark}
\begin{em}
By the Filippov Theorem, see e.g. \cite{AuFr}, Theorem 8.2.10, p. 316,
if $U$ is a complete metric space equipped with its Borel $\sigma$-algebra,
$g$ is measurable, $h$ is measurable bounded, with
$u\mapsto h(t,x,u)$  continuous on $U$, and if
$\Gamma$ takes  non-empty values, then $\Gamma$ admits a  measurable selection,
i.e. there exists a Borel measurable map $\Gamma
_{0}:[0,T]\times\bfC\times\R^d\to U$ such that
$\Gamma
_{0}\left(t,x,z\right)  \in\Gamma\left(t,x,z\right)  $ for $t\in [0,T]$, $x\in E,$ $z\in\R^d$.
\end{em}
\end{remark}

We are now ready to formulate the
assumptions we need .

\begin{hypothesis}\label{ipcosto}
\begin{enumerate}
 \item $(U,\calu)$ is a measurable space,  $g:U\to [0,\infty)$ is measurable,
 $h: [0,T]\times \bfC\times U\to\R^d$ is
measurable and bounded;
\item the hamiltonian $\psi$  defined in (\ref{hamiltoniana})
satisfies the requirements of  points 2 and 3 of  Hypothesis \ref{ip su psi};
\item the function $\phi:\bfC\rightarrow\R$ satisfies the requirements
of point 4 in Hypothesis
\ref{ip su psi}, namely it belongs to
$ \calg^{1}\left(  \bfC,\R\right)  $ and there exist $K>0$ and $m\geq0$ such that
\[
\left\vert \nabla\phi(x)  \right\vert \leq
K\left(  1+\left\vert x\right\vert _{\bfC}\right)^{m},\qquad x\in \bfC.
\]
\end{enumerate}
\end{hypothesis}

\begin{remark}
 \begin{em} \begin{enumerate}
 \item Hypothesis \ref{ipcosto} is stronger than  Hypothesis \ref{ip su psi}.
Indeed, point 1 of Hypothesis \ref{ip su psi} is a
 straightforward consequence of the fact that  $h$ is assumed to be bounded.

 \item
In the case $U\subset \R^k$, $h(t,x,u)=u$, the previous assumptions
require in particular that the set  $U$  where control processes take values should be  bounded.

\item The assumptions on the hamiltonian function $\psi$ can be easily
verified in specific cases. For instance if
$U$ is a closed ball of  $\R^k$ centered at the origin, and $g(u)=g_0(\vert u\vert^p)$
for some $p>1$ and some convex function $g_0:[0,\infty)\to [0,\infty)$ such that
$g\in C^{1}([0,\infty))$  and $g^{\prime}\left(
0\right)  >0$,
 then the hamiltonian is differentiable with respect to $z$  and
 $\psi$ satisfies points 2 and 3 of Hypothesis \ref{ip su psi}.
\end{enumerate}
\end{em}
\end{remark}

Now let us consider a probability space
$(  \widetilde{\Omega},\widetilde{\calf},\widetilde{\mathbb{P}})  $,
a standard
Wiener process $  \widetilde{W}  $
in $\R^d$, and
 the following forward-backward system:
\begin{equation}\label{fbsdetilde}
\left\{\begin{array}{l}\dis dy_\tau =
b(\tau,{X}_\tau)\; d\tau
+\sigma(\tau,{X}_\tau)\; d\widetilde{W}_\tau,\quad \tau\in
[t,T]\subset [0,T],
\\\dis
X_t=x,
\\\dis
 dY_\tau=\psi(X_\tau,Z_\tau)\;d\tau+Z_\tau\;d\widetilde{W}_\tau,
  \\\dis
  Y_T=\phi(X_T).
\end{array}\right.
\end{equation}
By Remark \ref{remark v}, the function
$v:\left[  0,T\right]  \times \bfC\rightarrow\mathbb{R}$
defined by the equality
\begin{equation}
v\left(
t,x\right)  =Y^{t,x}_t  \label{def v}
\end{equation}
is a functional of the coefficients $b,\sigma,\psi,\phi$ and
does not depend on the particular choice of
$(  \widetilde{\Omega},\widetilde{\calf},\widetilde{\mathbb{P}})  $
nor on the Wiener process
$  \widetilde{W} $.

In the following proposition we show  that the function $v$, defined in this way
by means of an appropriate forward-backward stochastic differential system,
plays a basic role in the control problem.

\begin{proposition}
\label{prop rel fond}Assume that Hypotheses \ref{ipotesi1} and \ref{ipcosto} hold true, and
that the cost functional is given in (\ref{funzcostofinale}).
 Let $v$ be defined in (\ref{def v}). Then for
every $t\in\left[  0,T\right]  $ and $x\in \bfC$ and for every admissible control
$u(\cdot)$ we have $J\left(  t,x,u(\cdot) \right)  \geq v\left(  t,x\right)  $.
\end{proposition}

\noindent {\bf Proof.} We fix $t,x$ and a control $u(\cdot)$.
Let $(X_{\tau}^{u})_{\tau \in [t,T]}$ be the corresponding process defined by
  (\ref{segmentou}). We define the process
\[
W_{\tau}^{u}=W_{\tau}+\int_{t\wedge\tau}^{\tau}h\left( s,X_s
^{u},u_s\right)  ds,\text{ \ \ }\tau\in\left[  0,T\right]  ,
\]
and we note that $X^{u}$ solves the equation
\begin{equation}\label{forward eq fond}
\left\{
\begin{array}
[c]{l}
dy^u_\tau =
b(\tau,X^u_{\tau})\; ds
+\sigma(\tau,X^u_{\tau})\;
 dW^u_\tau
,\text{ \ \ \ }\tau\in\left[  t,T\right] ,\\
X^u_t(\theta)=x(\theta),\quad \theta\in [-r,0].
\end{array}
\right.
\end{equation}
Since $h$ is bounded, we can apply the Girsanov theorem and deduce that there
exists a probability measure $\mathbb{P}^{u}$ on $\left(  \Omega
,\mathcal{F}\right)  $ such that $W^{u}$ is a Wiener process with respect to
$\mathbb{P}^{u}$.
We remark that, by uniqueness, $X^u$ is in fact a continuous process adapted to the
natural filtration generated by $W^u$ and augmented by the $\P^u$-null sets.
In $\left(  \Omega,\mathcal{F},\mathbb{P}^{u}\right)  $\ let
us consider the backward equation for the unknown process $\left(  Y_{\tau
}^{u},Z_{\tau}^{u}\right)  $, $\tau\in\left[  t,T\right]  $:
\begin{equation}
Y_{\tau}^{u}+\int_{\tau}^{T}Z_s^{u}dW_s^{u}=\phi\left(
X_{T}^{u}\right)+\int_{\tau}^{T}\psi\left(s,X_s^{u},Z_s
^{u}\right)  ds,\ \ \ \tau\in\left[  t,T\right]  .
\label{backward eq fond}
\end{equation}
We notice that the forward-backward system (\ref{forward eq fond})-(\ref{backward eq fond})
has the form (\ref{fbsdetilde}) and we conclude that
 $Y^{u}_t  =v\left(  t,x\right)  $, where $v$
is defined in (\ref{def v}); in particular, it does not depend on $u(\cdot)$.

Now we wish to prove that $\int_{t}^{T}
Z_s^{u}dW_s$ has finite expectation, equal to zero.
By the Burkhol\-der-Davis-Gundy inequalities, it is enough to
 prove
that
\begin{equation}\label{claimbdg}
    \mathbb{E}\left(  \int_{t}^{T}\left\vert Z_s^{u}\right\vert
^{2}ds\right)  ^{1/2}<\infty.
\end{equation}
We remember that
\[
\frac{d\mathbb{P}^{u}}{d\mathbb{P}}=\exp\left(  -\int_{t}^{T}h\left(s,
X_s^{u},u_s\right)  dW_s-\frac{1}{2}\int_{t}
^{T}\left\vert h\left(s,X_s^{u},u_s\right)  \right\vert
^{2}ds\right)  .
\]
We denote $\dfrac{d\mathbb{P}^{u}}{d\mathbb{P}}$ by $\rho$,
and by $\E^u$ the expectation with respect to $\P^u$. We estimate
\begin{align*}
\mathbb{E}\left(  \int_{t}^{T}\left\vert Z_s^{u}\right\vert
^{2}ds\right)  ^{1/2}  &  =\mathbb{E}^{u}\left[  \left(  \int_{t}
^{T}\left\vert Z_s^{u}\right\vert ^{2}ds\right)  ^{1/2}
\rho^{-1}\right] \\
&  \leq\left( \mathbb{E}^{u} \int_{t}^{T}\left\vert Z_s
^{u}\right\vert^{2}ds\right)  ^{1/2}\left(\mathbb{E}^{u}\left[
\rho^{-2}\right]  \right)^{1/2}.
\end{align*}
Since the process $Z^u$, solution to (\ref{backward eq fond}), is square-summable
withe respect to $\P^u$, it
 remains to prove that $\mathbb{E}^{u}\left[  \rho^{-2}\right]  $ is finite.
Noting that
\[
  \rho^{-1}  =  \exp\left(
\int_{t}^{T}h\left( s,X_s^{u},u_s\right)  dW_s
^{u}-\frac{1}{2}\int_{t}^{T}\left\vert h\left(s,X_s
^{u},u_s\right)  \right\vert^{2}ds\right) ,
\]
and recalling that $h$ is bounded
we get, for some constant $C$,
\begin{align*}
\mathbb{E}^{u}\left[  \rho^{-2}\right]   &  =\mathbb{E}^{u}\left[
\exp2\left(  \int_{t}^{T}h\left( s,X_s^{u},u_s\right)
dW_s^{u}-\frac{1}{2}\int_{t}^{T}\left\vert h\left(s,X_s^{u},
u_s\right)  \right\vert^{2}ds\right)  \right] \\
&  =\mathbb{E}^{u}\left[  \exp\left(  \int_{t}^{T}2h\left(s,X_s^{u},
u_s\right)  dW_s^{u}-\frac{1}{2}\int_{t}^{T}4\left\vert
h\left(s,X_s^{u},u_s\right)  \right\vert
^{2}ds\right)  \right. \\
&  \left.  \exp\left(  \int_{t}^{T}2\left\vert h\left(s,X_s
^{u},u_s\right)  \right\vert^{2}ds\right)  \right] \\
&  \le C\;\mathbb{E}^{u}\left[  \exp\left(  \int_{t}^{T}2h\left(s,X_s^{u},
u_s\right)  dW_s^{u}-\frac{1}{2}\int_{t}^{T}4\left\vert
h\left(s,X_s^{u},u_s\right)  \right\vert
^{2}ds\right)  \right. \\
&= C .
\end{align*}
Now (\ref{claimbdg}) is proved and therefore
$\int_{t}^{\tau}Z^u_sdW_s$  has zero expectation with respect to the
original probability $\mathbb{P}$. If  we set
$\tau=t$ in (\ref{backward eq fond}) and we take expectation with respect to
$\mathbb{P}$, we obtain
\begin{equation*}
v\left(  t,x\right) =\E\phi\left(  X_{T}^{u}\right)  +\mathbb{E}
\int_{t}^{T}\left[  \psi\left(s,X_s^{u},Z_s^{u}\right)
-Z_s^{u}h\left( s,X_s^{u},u_s\right)  \right]
ds.
\end{equation*}
Adding and subtracting $\E\int_{t}^{T}
g\left( u_s\right)  ds$ we arrive at
\begin{equation}
v\left(  t,x\right)  =J\left(t,x,u(\cdot)\right)
+\mathbb{E}\int_{t}^{T}\left[  \psi\left( s,X_s^{u},
Z_s^{u}\right)-Z_s^{u}h\left(s,X_s^{u},u_s\right)
-g( u_s)\right]  ds.
\label{backward rel fond}
\end{equation}
By the definition of $\psi$ the term in the square brackets is non positive
and consequently
$v\left(  t,x\right)  \le
J\left(  t,x,u (\cdot) \right) $.
\qed

The equality (\ref{backward rel fond}) can be regarded as a version of the
so-called {\em fundamental relation}. We immediately deduce the
following consequences:
\begin{proposition}
\label{propcontrolloforte}Let $t\in\left[  0,T\right]  $ and $x\in \bfC$\ be
fixed. If a control $u(\cdot)$ satisfies  $J\left( t,x,u(\cdot) \right)  =v\left(  t,x\right)  $
then $u(\cdot)$ is optimal for the control problem starting
from $x$ at time $t$.

Assume that the set-valued map $\Gamma$ has non empty values and it
admits a measurable selection $\Gamma
_{0}:[0,T]\times\bfC\times\R^d\to U$, and assume that
a control $u(\cdot)$ satisfies
\begin{equation}
u_{\tau}=\Gamma_{0}\left(\tau,X_{\tau}^{u},Z_{\tau}^{u}\right),  \text{
\ \ \ \ }\mathbb{P}\text{-a.s. for almost every }\tau\in\left[  t,T\right]  .
\label{u in gammazero}
\end{equation}
Then $J\left(t,x,u(\cdot)\right)  =v\left(  t,x\right)  $,
$u(\cdot)$ is optimal, and the optimal pair $(u(\cdot),X^u)$ satisfies the feedback law
\begin{equation}
u_{\tau}=\Gamma_{0}\left(\tau,X_{\tau}^{u},\nabla_0 v(\tau,X^u_\tau )\sigma(\tau,X^u_\tau)\right),
  \text{
\ \ \ \ }\mathbb{P}\text{-a.s. for almost every }\tau\in\left[  t,T\right]  .
\label{u-feedback}
\end{equation}
\end{proposition}

We note that
(\ref{u-feedback})
follows from  (\ref{u in gammazero}) and
(\ref{identificazioneZ}).

However, we can not
prove the existence of an optimal control satisfying (\ref{u in gammazero}) (and hence
(\ref{u-feedback})). Such a control can be shown to exist if there exists
a solution to the so-called closed-loop equation
\begin{equation}
\left\{\!\!
\begin{array}
[c]{l}
dy_{\tau}=b\left(  \tau,X_{\tau}\right)
d\tau+\sigma(  \tau,X_{\tau})[ h(\tau,X_\tau, \Gamma_{0}\left(\tau,X_{\tau},\nabla_0
v\left(\tau,X_{\tau}\right)\sigma \left(\tau,X_{\tau}\right)\right))d\tau+dW_{\tau}],
\;\tau\in\left[  t,T\right],  \\
X_{t}(\theta)=x(\theta),\qquad \theta \in [-r,0],
\end{array}
\right.  \label{closed loop eq}
\end{equation}
since in this case one can define an optimal control setting
$$
u_{\tau}=\Gamma_{0}\left(\tau,X_{\tau},\nabla_0 v(\tau,X_\tau )\sigma(\tau,X_\tau)\right).
$$
However, under the present assumptions, we can not
guarantee that the closed-loop equation has
a solution in the usual strong sense. To circumvent this difficulty
we will revert to  a weak formulation of the optimal control problem.

\subsection{Weak formulation of the optimal control problem}

We formulate the optimal control problem in the weak sense following the
approach of \cite{FlSo}, see e.g. chapter III. The main advantage
is that we will be able to solve the
closed loop equation in a weak sense, and hence to find an optimal control,
even if the feedback law is non smooth.

Initially, we are given the set $U$ and the functions $b,\sigma, h, g, \phi$.
By an {\em admissible control system}
we mean
$$(\Omega,\mathcal{F},
\left(\mathcal{F}_{t}\right)  _{t\geq 0}, \mathbb{P}, W,
u(\cdot),X^u),$$
where
$(  \Omega,\mathcal{F},\left(  \mathcal{F}_{t}\right)  _{t\geq0},\mathbb{P})  $
is a filtered probability space satisfying the usual conditions,  $W$ is
 an $\R^d$-valued standard Wiener process with respect to
 $(\mathcal{F}_{t})$ and $\mathbb{P}$,
 $u$ is an $(  \mathcal{F}_{t})  $-predictable process with values in $U$,
 $X^u$ satisfies
(\ref{segmentou})-(\ref{SDDEcontrollata1}).
An admissible control system will be  briefly
denoted by $\left(W,u,X^{u}\right)$ in the following.
 Our aim is now to minimize the cost functional
\begin{equation}\label{funzcostofinaledebole}
J\left(  t,x,\left(  W,u,X^{u}\right)  \right)
=\E\int_t^Tg(u_s)\; ds
+
\E\,\phi(X^u_T)
\end{equation}
over all the admissible control systems $\left(W,u,X^{u}\right)$. We can prove the following results:
\begin{theorem}
\label{teo su controllo debole}
Assume that Hypoteses \ref{ipotesi1} and \ref{ipcosto} hold true, and
that the cost functional is given in (\ref{funzcostofinaledebole}).
 Let $v$ be defined in (\ref{def v}). Then for
every $t\in\left[  0,T\right]  $ and $x\in \bfC$
and for all admissible control system $\left(
W,u,X^{u}\right)  $ we have
\[
J\left(  t,x,\left(  W,u,X^{u}\right)  \right)  \geq v\left(  t,x\right)
\]
and the equality holds if and only if
\begin{equation*}
u_{\tau}\in\Gamma(\tau,X_{\tau}^{u},\nabla_0 v(\tau,X^u_\tau )
\sigma(\tau,X^u_\tau)),\text{ \ }\mathbb{P}\text{-a.s. for a.a. }\tau
\in\left[  t,T\right]  .
\end{equation*}

Moreover
assume that the set-valued map $\Gamma$ has non empty values and it
admits a measurable selection $\Gamma
_{0}:[0,T]\times\bfC\times\R^d\to U$.
Then an admissible control system $\left(
W,u,X^{u}\right)  $  satisfying the feedback law
\begin{equation*}
u_{\tau}=\Gamma_{0}(X_{\tau}^{u},\nabla_0 v(\tau,X^u_\tau )
\sigma(\tau,X^u_\tau)),\text{ \ }\mathbb{P}\text{-a.s. for a.a. }\tau
\in\left[  t,T\right]
\end{equation*}
is optimal.

Finally, the closed loop equation (\ref{closed loop eq})
admits a weak solution
$(\Omega,\mathcal{F},
\left(\mathcal{F}_{t}\right)  _{t\geq 0}, \mathbb{P}, W,X)$
which is unique in law and setting
$$
u_{\tau}=\Gamma_{0}\left(\tau,X_{\tau},\nabla_0 v(\tau,X_\tau )\sigma(\tau,X_\tau)\right),
$$
we obtain an optimal admissible control system $\left(
W,u,X\right)  $.
\end{theorem}

\noindent {\bf Proof.} The proof follows from the fundamental relation
(\ref{backward rel fond}) and the same arguments
leading to Proposition \ref{propcontrolloforte} and the remarks following it.
The only difference here is the solvability
of the closed loop equation in a weak sense, which is however a
standard application of  a Girsanov change of
measure.
\qed

\section{Parabolic equations and application to a pricing problem}
\label{sez-nlk}
Let us consider again the Markov process $\{X^{t,x}_\tau,\, 0\le t\le \tau\le T,x\in\bfC\}$,
defined by the formula
(\ref{defdixdue}), starting from the
family of solutions
to equation
(\ref{eqscalare-s}). Let us denote by
$(\call_t)_{t\in[0,T]}$ the corresponding generator. Thus, each
$\call_t$ is a second order differential
operator acting on a suitable domain consisting of real functions
defined on $\bfC$.
  In the autonomous case, a description of the generator, denoted
by $\call$, was given
in Section 2, remark \ref{remarkL}.
In this section we treat semilinear parabolic equations driven by
$(\call_t)$, which are generalizations of the Kolmogorov equations.
We will introduce a concept of solution, called
mild solution, that does not require a description of the generators.
In the sequel the notation $\call_t$ will be used  only in a formal way.

The parabolic equations that we study have the following form:
\begin{equation}\label{nlkdeg}
  \left\{\begin{array}{l}\dis
\frac{\partial v(t,x)}{\partial t}+\call_t v(t,x) =
\psi (t,x,v(t,x),\nabla_0 v(t,x)\,\sigma(t ,x)),\\ \\
\dis v(T,x)=\phi(x), \qquad t\in [0,T],\;x\in \bfC,
\end{array}\right.
\end{equation}
with unknown function $v:[0,T]\times \bfC\to \R$ and given coefficients
$\psi:[0,T]\times\bfC\times\R\times\R^d
\rightarrow\mathbb{R}$ and $\mathbb{\phi}:\bfC\rightarrow\mathbb{R}$.
We recall the notation $\nabla_0 v(t,x)$ introduced in
(\ref{defnablazero}).
We note in particular that $\nabla_0v(t,x)$ is a vector in $\R^n$ whose components
are denoted
$\nabla_0^kv(t,x)$ ($k=1,\ldots,n$). If we denote
 by
 $\sigma^i_k(t,x)$ ($k=1,\ldots, n$, $i=1,\ldots,d$) the components of the matrix
$\sigma(t,x)$,
then $\nabla_0 v(t,x)\,\sigma(t ,x)$ denotes the vector in $\R^d$ whose
components are $\sum_{k=1}^n \nabla_0^kv(t,x)\sigma^i_k(t,x)$,
($i=1,\ldots,d$).

Recalling the definition of the transition semigroup $P_{t,\tau}$ given in
(\ref{defdiP}), and writing the variation of constants formula for a solution
to  (\ref{nlkdeg}),
we formally obtain
\begin{equation}\label{mild sol hjb}
 v(t,x) =P_{t,T}[ \phi](x)
  -\int_t^TP_{t,\tau }\Big[
\psi ( \cdot,v(\tau ,\cdot), \nabla_0v(\tau ,\cdot)\,
\sigma(\tau ,\cdot))\Big](x) \;
d\tau \qquad t\in [0,T],\,x\in\bfC.
\end{equation}
We notice that this formula is meaningful if  $\nabla_0 v $ is well defined
and provided $\phi$ and $\psi$
satisfy some growth and measurability conditions. This way we arrive at
the
following definition of mild solution of the semilinear Kolmogorov equation
(\ref{nlkdeg}).

\begin{definition}
\label{defsolmildkolmo}A function $v:\left[  0,T\right]  \times \bfC\rightarrow
\mathbb{R}$ is a mild solution of the semilinear Kolmogorov equation
(\ref{nlkdeg}) if
$v \in \calg^{1}\left([0,T]\times  \bf C;\mathbb{R}
\right)  $, there exist $C>0, q\geq 0$ such that
\begin{equation}\label{iposuv}
   |v(t,x)|+|\nabla_x v(t,x)|\leq C(1+|x|)^{q},
   \qquad   t\in [0,T], \; x\in \bfC,
\end{equation}
and the equality (\ref{mild sol hjb}) holds.
\end{definition}

The space
$\calg^{1}\left([0,T]\times  \bf C;\mathbb{R}
\right)  $ was described in Remark (\ref{spaziosol}). $|\nabla_x v(t,x)|$
denotes the total variation norm of the $\R^n$-valued finite Borel measure
$\nabla_x v(t,x)$ on $[-r,0]$.

\begin{theorem}\label{teo-nlkdeg} Assume that Hypotheses
\ref{ipotesi1} and \ref{ip su psi} hold true. Then there exists a unique mild solution $v$ of
(\ref{nlkdeg}). The function $v$ coincides with the one introduced in
Corollary \ref{proprietadiv}.
\end{theorem}
\noindent{\bf Proof.} At first we prove existence. For fixed
$t\in \left[  0,T\right]$ and $x\in \bfC$, let
$(X_\tau^{t,x}, Y_\tau^{t,x}, Z_\tau^{t,x})_{\tau\in[0,T]}$ denote the
solution of the forward-backward system (\ref{fbsde}) and
let $v(t,x)$ be defined by equality (\ref{defdiv}).
The required regularity and growth conditions of the function $v$ were
proved in
Corollary \ref{proprietadiv}, so it remains
to prove that $v$ satisfies equality (\ref{mild sol hjb}). To this
aim we evaluate
\begin{align*}
P_{t,\tau}\left[  \psi(\cdot,v\left(  \tau,\cdot\right),
\nabla_0 v\left(  \tau,\cdot\right)\sigma \left(  \tau,\cdot\right))\right]
\left(  x\right)
&  =\mathbb{E}\left[  \psi\left(\tau, X_\tau^{t,x}, v(\tau, X_\tau^{t,x} ),
\nabla_0 v(\tau,Y_\tau^{t,x})\sigma(\tau,X_\tau^{t,x}) \right) \right] \\ \nonumber
&  =\mathbb{E}\left[  \psi\left(X_\tau^{t,x}, Y_\tau^{t,x},
Z_\tau^{t,x}  \right)  \right],
\end{align*}
where the last equality follows from
(\ref{identificazioneY}) and
(\ref{identificazioneZ}).
In
particular we obtain
\begin{equation}
\int_{t}^{T}P_{t,\tau}\left[  \psi(\cdot,v\left(  \tau
,\cdot\right),\nabla_0 v\left(  \tau,\cdot\right)\sigma(\tau,\cdot)\right]
\left(  x\right)  d\tau=\mathbb{E}\int_{t}^{T}%
\psi\left(  \tau,X_\tau^{t,x},Y_\tau^{t,x}  ,Z_\tau^{t,x}\right)
d\tau. \label{psi semigruppo}
\end{equation}
Since the pair $\left(  Y^{t,x},Z^{t,x}\right)  $ is a solution to the backward equation
(\ref{bsde}) we have
\[
Y_t^{t,x}  +\int_{t}^{T}Z_\tau^{t,x}  dW_{\tau
}=\phi\left(  X_T^{t,x}\right)  +\int_{t}^{T}
\psi\left(X_\tau^{t,x},Y_\tau^{t,x} ,Z_\tau^{t,x} \right)
d\tau.
\]
Taking expectation and applying formula (\ref{psi semigruppo}) we get the
equality (\ref{mild sol hjb}).

It remains to prove uniqueness. Let $v$ be a mild solution
 to (\ref{nlkdeg}), so that  for every $s\in\left[  t,T\right] \subset \left[  0,T\right] $,
\[
v\left(  s,x\right)  =\mathbb{E\phi}\left(  X_T^{s,x}\right)
+\mathbb{E}\int_{s}^{T}\psi\left(  X_\tau^{s,x},v(\tau,X_\tau^{s,x})  ,\nabla_0
v\left(  \tau,X_\tau^{s,x}\right)\sigma(\tau,X_\tau^{s,x})\right) d\tau.
\]
Since $X_\tau^{s,x} $ is independent on $\mathcal{F}_{s}$,
the expectations occurring in the last formula can be replaced
by conditional expectations given
$\mathcal{F}_{s}$. Next we note that $x$ can be replaced by $X_s^{t,x} $,
since $X_s^{t,x}$ is $\mathcal{F}_{s}$-measurable. Using the identity
$X_\tau^{s,X_s^{t,x}}=X_\tau^{t,x} $, which follows easily from uniqueness
of the forward equation,
we finally obtain
\begin{align*}
v\left(  s,X_s^{t,x}\right)   &  =\E^{\mathcal{F}_{s}%
}\mathbb{\phi}\left(  X_T^{t,x}  \right)  +
\E^{\mathcal{F}_{s}}\int_{s}^{T}\psi\left(X_\tau^{t,x},v(\tau,X_\tau^{t,x})
,\nabla_0v\left(  \tau,X_\tau^{t,x}\right)\sigma\left(  \tau,X_\tau^{t,x}\right)\right)  d\tau\\
&  =\mathbb{E}^{\mathcal{F}_{s}}\mathbb{\eta}-\mathbb{E}^{\mathcal{F}_{s}}%
\int_{t}^{s}\psi\left(X_\tau^{t,x},v(\tau,X_\tau^{t,x})  ,\nabla_0v\left(
\tau,X_\tau^{t,x}\right)\sigma\left(  \tau,X_\tau^{t,x}\right)\right)  d\tau,
\end{align*}
where we have defined $\eta=\phi\left(  X_T^{t,x}\right)  +
\dis\int_{t}^{T}
\psi\left( X_\tau^{t,x},v(\tau,X_\tau^{t,x}) ,\nabla_0
v\left(  \tau,X_\tau^{t,x}  \right)\sigma(\tau,X_\tau^{t,x})\right)d\tau$. By the representation theorem of
martingales, see e.g. \cite{DP1}, theorem 8.2, there exists
a predictable process $\widetilde{Z}\in L^{2}
\left(  \Omega\times\left[  0,T\right]  ,\R^d
\right)  $, such that $\E^{\mathcal{F}_{s}}\eta=
\dis\int_{t}^{s}
\widetilde{Z}_{\tau}dW_{\tau}+v\left(  t,x\right) $, $s\in [t,T]$. So
\begin{equation}
v\left(  s,X_s^{t,x} \right)  =v\left(  t,x\right)  +\int_{t}
^{s}\widetilde{Z}_{\tau}dW_{\tau}-\int_{t}^{s}\psi\left(X_\tau^{t,x},
v(\tau,X_\tau^{t,x})  ,\nabla_0v\left(  \tau,X_ \tau^{t,x}\right)
\sigma(\tau,X_\tau^{t,x})  \right) d\tau. \label{v semimartingala}
\end{equation}
Now we compute the joint quadratic variation with $W^i$ of the processes
occurring at both sides of this equality,
on an interval $[t,T']\subset [t,T)$.
Considering the right-hand side we obtain
$\int_{t}
^{T'}\widetilde{Z}^i_{\tau}d{\tau}$
by the rules of stochastic calculus.
By Theorem \ref{esistejointvar} we have
$
\<v(\cdot,X_\cdot^{t,x}),W^i\>_{[t,T']} =\int_t^{T'}
\sigma^i(\tau,X_\tau^{t,x})\nabla_0 v(\tau,X_\tau^{t,x})\;d\tau
$.
Therefore we have $\widetilde{Z}_{\tau}=\sigma(\tau,X_\tau^{t,x})\nabla_0 v(\tau,X_\tau^{t,x})$
and
equality (\ref{v semimartingala}) can be rewritten as
$$
\begin{array}{lll}
  v(  s,X_s^{t,x})
&  =&\dis v\left(  t,x\right)  +\int_{t}^{s}\nabla_{0}v\left(  \tau,X_\tau^{t,x}
\right)\sigma(\tau,X_\tau^{t,x})  dW_{\tau}
\\&&\dis
-\int_{t}^{s}\psi\left( X_\tau^{t,x},
v(\tau,X_\tau^{t,x}),\nabla_{0}v\left(\tau,X_\tau^{t,x}\right)
\sigma(\tau,X_\tau^{t,x})\right)  d\tau\\
&  =&\dis \phi\left( X_T^{t,x}\right)  -\int_{s}^{T}\nabla
_0v\left(  \tau,X_\tau^{t,x} \right)\sigma(\tau,X_\tau^{t,x})dW_{\tau}
\\
&&\dis
+\int_{s}^{T}\psi\left(X_\tau^{t,x},v(\tau,X_\tau^{t,x}),\nabla_0v\left(
\tau,X_\tau^{t,x}  \right)\sigma(\tau,X_\tau^{t,x})\right)  d\tau.
\end{array}
$$
By comparing with the backward equation in (\ref{fbsde}) we see that the
pairs of processes
$(  Y_ s^{t,x}  ,Z_ s^{t,x})$
and
$
\left(  v(  s,X_s^{t,x}  )  ,\nabla_0v(
s,X_s^{t,x} )\sigma(s,X_s^{t,x})\right)
$, $  {s\in [t,T]}$,
solve the same equation. By uniqueness of the solution
we have $Y_s^{t,x} =v(
s,X_s^{t,x})  ,$ $s\in [t,T]$, and for $s=t$ we get  $Y_t^{t,x}=v\left(  t,x\right)  $.
\qed

\begin{remark}{\em
The proof of uniqueness is based on an application of Theorem \ref{esistejointvar}.
Inspection of the proof shows that uniqueness holds in a larger class of functions.
Namely, if a
Borel measurable functions
$v:[0,T]\times \bfC\to \R$ satisfies
 $v(t,\cdot)\in \calg^1(\bfC,\R)$ for every $t\in[0,T]$,
and the inequality
$$
   |v(t,x)|+|\nabla_x v(t,x)|\leq C(1+|x|)^{q},
   \qquad    t\in [0,T],x\in E,
$$
holds
for some $C>0, q\geq 0$, and the equation  (\ref{mild sol hjb})
holds, then $v$ coincides with the solution constructed
in Theorem \ref{esistejointvar}.
}
\end{remark}

\begin{remark}
\emph{ If in (\ref{nlkdeg}) we take as $\psi$ the hamiltonian
defined in (\ref{hamiltoniana}) and as $\phi$ the final
cost in the cost functional (\ref{funzcostofinale}), equation (\ref{nlkdeg}) is
the Hamilton Jacobi Bellman equation related to the stochastic optimal control problem we
have treated in section 5.
It turns out that the value function coincides with the solution of (\ref{nlkdeg})}.
\end{remark}

\subsection{Application to pricing}
We consider a financial  market, of Black and Scholes type, with one risky
asset, whose price at time $t$ is denoted by $S_t$,
and one non risky asset, whose price is
denoted by $B_t$. We assume the following prices evolution:
\begin{equation}\label{black-sholes}\left\{\begin{array}{ll}
\dis d S_t = \mu(t,S_{t+\cdot})\;S_t\; dt \; + \;
  \sigma(t,S_{t+\cdot})\;S_t\; dW_t,& t\in [0,T],\;
  \\ S_\theta=s_\theta,&\theta\in
[-r,0],
 \\
   \dis d B_t = \rho B_t\; dt,&t\in [0,T], \\ B_0=1, &
\end{array}\right.
\end{equation}
where $\rho>0$, $r>0$ and $s\in
\bfC= C([-r,0],\R)$.
We notice that the coefficients $\mu$ and $\sigma$ depend on the
past trajectory: $S_{t+\cdot}$ stands for the past trajectory of length  $r$,
i.e. $S_{t+\cdot}=(S_{t+\theta})_{\theta\in[-r,0]}$.
Moreover we consider a contingent claim of the form
$$\phi(S_{T+\cdot}),$$
where $\phi:\bfC\to \R$.
If $r>T$ then the claim
depends on the whole evolution in time of the prices
of the shares: see \cite{bj}, \cite{MuRu} or \cite{will} and
references within for a general discussion on such kind of options,
usually referred to as {\em path-dependent}.

We denote by $\pi_t$ the value of the investor's portfolio
invested in the risky asset at time $t$. $\pi$ is called a trading strategy;
we will only consider
predictable trading strategies which are square-integrable, i.e.
$\E\int_0^T|\pi_t|^2dt<\infty$.
We notice that the value $V_t$
of the corresponding self-financing portfolio satisfies the
equation
\begin{equation}\label{portfolio-eq}
 d V_t = \rho V_t\;dt +\pi_t\; \sigma(t,S_{t+\cdot})\;\theta (t,S_{t+\cdot})\; dt+
   \pi_t \;\sigma(t,S_{t+\cdot})\; dW_t,
\end{equation}
where
\[
\label{riskpremium}
\theta(t,S_{t+\cdot})=\frac{\mu(t,S_{t+\cdot})-\rho}{\sigma(t,S_{t+\cdot})}
\]
is called the risk premium.

At time $T$ the investor has to
pay a contingent claim of the form $\phi(S_{T+\cdot})$, where
$\phi : {\bf C}\rightarrow
\R$ is some given function. The pricing problem is to
 find and characterize pairs $(\pi,V_0)$ consisting of
a  strategy  $\pi$ and an initial capital $V_0\in\R$ such that
$$
V_T=\phi(S(T+\cdot)).
$$
$\pi$ is then called a hedging strategy and
$V_0$ is called the
fair price of the claim at time $t=0$.

Throughout this section we assume the following.
\begin{hypothesis}\label{ipotesipricing}
\begin{enumerate}
\item $(W_t)_{t\geq 0}$ is a real Wiener process defined
in a complete probability space $(\Omega,\mathcal{F},\P)$ and
$(\mathcal{F}_t)_{t\geq 0}$ is the filtration generated by $W$
augmented with null  sets.

\item $\mu:[0,T]\times  \bf C\to
\R$ is Borel measurable and bounded, and there exists $L>0$ such that
\begin{equation}\label{ipotesisumu}
| \mu(t,f^1)f^1(0)-\mu(t,f^2)f^2(0)|
\leq L \vert f_1-f_2\vert_{\bf C},
\end{equation}
for all $t\in [0,T]$, $f^1,f^2\in {\bf C}$; moreover,
$\mu(t,\cdot) \in\calg^1(\bf C,
\R)$ for all $t\in [0,T]$.

\item $\sigma :[0,T] \times {\bf C}\to
\R$ is Borel measurable and there exists $c>0$ such that
\begin{equation}\label{ipotesisusigmalimitatodasotto}
 |\sigma(t,f)|\geq c,
\end{equation}
for every $f\in{\bf C}$, so that the risk premium in (\ref{riskpremium})
is well defined and bounded; moreover
\begin{equation}\label{ipotesisusigma}
| \sigma(t,f^1)f^1(0)-\sigma(t,f^2)f^2(0)|
\leq L \vert f_1-f_2\vert_{\bf C},
\end{equation}
 for a suitable $L>0$ and for all
    $t\in [0,T]$, $f^1,f^2\in {\bf C}$; finally,
$\sigma(t,\cdot) \in\calg^1(\bf C,
\R)$ for all    $t\in [0,T]$

    \item $\phi\in \calg^1(\bfC,\R)$ satisfies $|\nabla\phi(x)|\le C (1+|x|_\bfC)^m$ for
    all $x\in\bfC$ and some $C>0$ and $m\ge 0$.
\end{enumerate}
\end{hypothesis}

By the Girsanov theorem there exists a probability measure, called risk-neutral
probability,
for which
$$
\overline{W}_t=\int_0^t \theta (\tau,S_{\tau+\cdot})\;d\tau + W_t, \qquad t\in [0,T],
$$
is a Wiener process. Then
$$
dS_t=\rho\,S_t\;dt+\sigma(t,S_{t+\cdot})\,S_t\;d\overline{W}_t,
\qquad
dV_t= \rho V_t\;dt+ \pi_t\sigma(t,S_{t+\cdot})\;d\overline{W}_t.
$$
The existence of a hedging strategy can be established as follows:
using the results of Section \ref{sezioneforback} we first
find a solution
 to the following
forward-backward stochastic differential system
\begin{equation}\label{fbsde_pricing}
    \left\{\begin{array}{l}\dis
    dS_t=\rho\;S_t\;dt+\sigma\,(t,S_{t+\cdot})\;S_t\;d\overline{W}_t,\quad t\in
[0,T],
\\\dis
S_{0+\cdot}=s,
\\\dis
dV_t= \rho V_t\;dt+ Z_t\;d\overline{W}_t,
  \\\dis
  V_T=\phi(S_{T+\cdot}).
\end{array}\right.
\end{equation}
Next, recalling
 (\ref{ipotesisusigmalimitatodasotto}),
we note that the required hedging strategy can be recovered from the process
$Z$ setting $\pi_t=Z_t/\sigma(t,S_{t+\cdot})$.

However, a better characterization of the hedging strategy  and the fair price
of the claim can be obtained. We first
consider,
for arbitrary $t\in [0,T]$ and $s\in\bfC$,  the  following
forward-backward system, which generalizes (\ref{fbsde_pricing}):
\begin{equation*}
    \left\{\begin{array}{l}\dis
    dS^{t,s}_\tau=\rho\,S^{t,s}_\tau\;dt+\sigma\,(\tau,S^{t,s}_{\tau+\cdot})
    \;S^{t,s}_\tau\;d\overline{W}_\tau,\quad \tau\in
[t,T],
\\\dis
S^{t,s}_{t+\cdot}=x,
\\\dis
dV^{t,s}_\tau= \rho V^{t,s}_\tau\;d\tau+ Z^{t,s}_\tau\;d\overline{W}_\tau,
  \\\dis
  V^{t,s}_T=\phi(S^{t,s}_{T+\cdot}),
\end{array}\right.
\end{equation*}
with unknown triple $(S_\tau^{t,s},V_\tau^{t,s}, Z_\tau^{t,s})$.
Setting  $X^{t,s}_\tau=S^{t,s}_{\tau+\cdot}$,
then $X$ is a Markov process in $\bfC$
with generator $\call$. We finally define
$$v(t,s)=V_t^{t,s},\qquad t\in [0,T],\;s\in\bfC.
$$
It follows from Corollary \ref{remark v}
that
$Z_\tau^{t,s}=\nabla_0 v(\tau,X_\tau^{t,x})\; \sigma (\tau,X_\tau^{t,x})$.
We conclude that the fair price and the hedging strategy are uniquely determined as
$$
V_0=v(0,s),\qquad \pi_t=\frac{Z_t^{0,s}}{\sigma(t, X_t^{0,s})}=\nabla_0 v(t,X_t^{0,s})=
{\nabla_0}
v(t,S_{t+\cdot}).
$$
Moreover, see Theorem \ref{teo-nlkdeg}, $v(t,s)$ is characterized as the unique mild solution of
the equation
\begin{equation}\label{kolmo_pricing}
  \left\{\begin{array}{l}\dis
\frac{\partial v(t,x)}{\partial t}+\call v(t,x) =
r\, v(t,x),\\ \\
\dis u(T,x)=\phi(x), \qquad t\in [0,T],\;
x\in \bfC,
\end{array}\right.
\end{equation}
which can be considered as a generalization of the Black-Scholes equation to the present setting.

\end{document}